\documentclass[reqno,12pt]{amsart}

\usepackage[all]{xy}
\usepackage{amsthm, amssymb, amsmath, amscd, mathrsfs}
\usepackage{geometry}
\usepackage[dvips]{graphicx, color}

\def\cal{\mathcal}

\def\NN{{\Bbb N}}
\def\PP{{\Bbb P}}

\def\ZZ{{\Bbb Z}}

\def\11{{1\kern-3.5pt 1}}
\def\mumu{{\mu\kern-4.2pt\mu}}

\def\boxtimes{\setbox0\hbox{$\Box$}\copy0\kern-\wd0\hbox{$\times$}}

\newcommand{\lotimes}{\otimes^{\bf{L}}}

\def\End{\operatorname {E nd}}
\def\Ext{\operatorname {Ext}}

\def\Hom{\operatorname {Hom}}
\def\id{\operatorname {id}}

\def\Ker{\operatorname {ker}}

\def\pd{{\operatorname {pd}}}

\def\Aut{\operatorname{Aut}}

\def\Bimod{\operatorname{Bimod}}

\def\Coker{\operatorname{Coker}}

\def\dim{\operatorname{dim}}

\def\End{\operatorname{End}}

\def\Ext{\operatorname{Ext}}

\def\uExt{\operatorname{\underline{Ext}}}

\def\gldim{\operatorname{gldim}}
\def\grmod{\operatorname{grmod}}

\def\GrMod{\operatorname{GrMod}}

\def\Hom{\operatorname{Hom}}
\def\RHom{\operatorname{RHom}}

\def\cHom{\operatorname{{\cH}om}}
\def\uHom{\operatorname{\underline{Hom}}}

\def\Im{\operatorname{Im}}

\def\Ker{\operatorname{Ker}}

\def\max{\operatorname{max}}

\def\mod{\operatorname{mod}}
\def\Mod{\operatorname{Mod}}

\def\Proj{\operatorname{Proj}}
\def\proj{\operatorname{proj}}

\def\sup{\operatorname{sup}}

\def\tails{\operatorname{tails}}
\def\Tails{\operatorname{Tails}}

\def\uTor{\operatorname{\underline{Tor}}}

\def\tors{\operatorname{tors}}
\def\Tors{\operatorname{Tors}}

\def\coh{\operatorname{coh}}

\def\uExt{\operatorname{\underline{Ext}}}

\def\uHom{\operatorname{\underline{Hom}}}
\def\RuHom{\operatorname{R\underline{Hom}}}
\def\uTor{\operatorname{\underline{Tor}}}

\def\R{\operatorname{R}}

\def\l{\leftarrow}

\let\oldtext\text
\def\text#1{\oldtext{\normalshape #1}}

\def\a{\alpha}
\def\b{\beta}
\def\c{\gamma}

\def\e{\epsilon}

\def\l{\lambda}
\def\s{\sigma}
\def\t{\tau}

\def\cA{{\cal A}}

\def\cC{{\cal C}}
\def\cD{{\cal D}}
\def\cE{{\cal E}}
\def\cF{{\cal F}}
\def\cG{{\cal G}}

\def\cHom{{\cal H}om}

\def\cL{{\cal L}}
\def\cM{{\cal M}}
\def\cN{{\cal N}}
\def\cO{{\cal O}}
\def\cP{{\cal P}}

\def\cX{{\cal X}}
\def\cY{{\cal Y}}

\def\Bimod{\operatorname{Bimod}}

\def\add{\operatorname{add}}

\def\<{\langle}
\def\>{\rangle}
\def\CM{\operatorname {CM}}

\def\sC{\mathscr C}
\def\sD{\mathscr D}
\def\sT{\mathscr T}

\def\uotimes{\underline{\otimes}}

\def\ulotimes{\underline{\otimes}^{\bf {\rm L}}}

\def\proj C{\add\{\cP_j\}_{j\in \ZZ}}

\newtheorem{lemma}{Lemma}[section]
\newtheorem{proposition}[lemma]{Proposition}
\newtheorem{theorem}[lemma]{Theorem}

\newtheorem{prop}[lemma]{Proposition}

{

}

\theoremstyle{definition}

\newtheorem{definition}[lemma]{\sl Definition}

\newtheorem{question}[lemma]{Question}

\theoremstyle{remark}

\newtheorem{remark}[lemma]{Remark}

\newcommand{\intotimesb}{{\underline{\otimes}}_{B}}

\newcommand{\intotimesc}{{\underline{\otimes}}_{C}}

\begin{document}

\pagenumbering{arabic}



\title{A categorical characterization of quantum projective $\mathbb{Z}$-spaces}

\author{Izuru Mori \and Adam Nyman}

\address{Department of Mathematics, Faculty of Science, Shizuoka University, Shizuoka 422-8529, JAPAN}

\email{mori.izuru@shizuoka.ac.jp}

\address{Department of Mathematics, 516 High Street, Western Washington University, 98225-9063, USA}

\email{adam.nyman@wwu.edu}



\thanks {The first author was supported by Grant-in-Aid for Scientific
Research (C) 20K03510 Japan Society for the Promotion of Science.}


\noindent


\begin{abstract}
In this paper we study a generalization of the notion of AS-regularity for connected $\ZZ$-algebras defined in \cite{MN}.  Our main result is a characterization of those categories equivalent to noncommutative projective schemes associated to right coherent regular $\ZZ$-algebras, which we call quantum projective $\ZZ$-spaces in this paper.

As an application, we show that smooth quadric hypersurfaces and the standard noncommutative smooth quadric surfaces studied in \cite{SV, MU1} have right noetherian AS-regular $\ZZ$-algebras as homogeneous coordinate algebras.  In particular, the latter are thus noncommutative $\PP^1\times \PP^1$ (in the sense of \cite{Vq}).
\end{abstract}

\maketitle

\tableofcontents

Throughout this paper, we work over a field $k$.

\section{Introduction}
In noncommutative algebraic geometry, one studies so-called noncommutative schemes from a geometric perspective. These schemes are often abelian categories with properties in common with categories of coherent sheaves over a scheme.  Those noncommutative schemes which behave like categories of coherent sheaves over {\it projective} schemes have particular significance, as they can be studied via their global invariants.  For this reason, it may be useful to characterize these abelian categories, and such a characterization constitutes the main result of this paper.

The first result along these lines was due to Artin and Zhang \cite[Theorem 4.5]{AZ}, who characterized those triples $(\sC, {\mathscr A}, s)$ consisting of a $k$-linear abelian category $\sC$, a distinguished object $\mathscr A$ (thought of as a structure sheaf), and an autoequivalence $s$ of the category $\sC$, which are equivalent to a triple of the form $(\tails A, \pi(A), s)$, where $A$ is some right noetherian $\mathbb{N}$-graded $k$-algebra satisfying a homological condition called $\chi_{1}$, $\grmod A$ is the category of finitely generated graded right $A$-modules, $\tors A$ is the full subcategory of $\grmod A$ consisting of right-bounded modules, $\tails A := \grmod A/\tors A$ is the quotient category with quotient functor $\pi$, and (abusing notation) $s$ is induced by the shift functor in $\grmod A$.  This result raises a natural question:

\begin{question} \label{quest.triple}
Is there a characterization of categories of the form $\tails A$ for suitably well behaved $\mathbb{Z}$-graded algebras $A$?
\end{question}

Artin and Zhang's result was later generalized by the first author and Ueyama \cite[Theorem 2.6]{MU1} to the case in which $A$ is right-coherent.  The first author and Ueyama then used this generalization to address Question \ref{quest.triple}.  In particular, they obtained a characterization of abelian categories equivalent to noncommutative projective schemes with homogenous coordinate ring a graded right coherent AS-regular algebra over a finite dimensional algebra $R$ of finite global dimension \cite[Theorem 4.1]{MU1}.

In a separate development,  Bondal and Polishchuk introduced the notion of $\mathbb{Z}$-algebra \cite{BP}, and illustrated the utility of this concept in the study of $\mathbb{Z}$-graded algebras.  Sierra provided further evidence \cite{S} that working with $\mathbb{Z}$-algebras simplifies aspects of the theory of $\mathbb{Z}$-graded algebras.  On the other hand, much of the theory of modules over $\mathbb{Z}$-graded algebras can be generalized to the $\mathbb{Z}$-algebra context (see, for example \cite{MN}).

Many noncommutative projective schemes with $\mathbb{Z}$-algebra coordinate rings have been studied. For example, Van den Bergh discovered notions of noncommutative $\PP^1\times \PP^1$
\cite{Vq} and noncommutative $\mathbb{P}^{1}$-bundles over a pair of smooth schemes \cite{Vn}.  Specializing the latter construction to the case where the base schemes are spectra of fields (or division rings), one obtains the notion of a noncommutative projective line.  Polishchuk \cite{P} found sufficient conditions for a $k$-linear abelian category to be of the form $\tails A$ for a right coherent positively graded $\mathbb{Z}$-algebra $A$.  He applied this result to construct $\mathbb{Z}$-algebra homogenous coordinate rings for noncommutative elliptic curves \cite{P2}.  Efimov, Lunts and Orlov constructed noncommutative Grassmannians with $\mathbb{Z}$-algebra homogeneous coordinate rings \cite{ELO}, providing further evidence for the significance of the notion of $\mathbb{Z}$-algebra.

Returning to characterizations of noncommutative projective  schemes, in \cite{Na}, the second author characterized those abelian categories equivalent to noncommutative projective lines over a pair of division rings.  The purpose of this paper is to obtain a $\mathbb{Z}$-algebra version of \cite[Theorem 4.1]{MU1} which generalizes \cite[Theorem 4.2]{Na}.

Instead of characterizing categories equivalent to $\tails A$ where $A$ is an AS-regular $\mathbb{Z}$-algebra, we use a related notion of regularity, called ASF$^{++}$-regularity (see Section \ref{section.regularity}).  In the $\mathbb{Z}$-graded case, AS-regularity implies that, if $\tau$ is the torsion functor, then $D\R\t(A)\cong A(-\ell)_{\nu}[d]$ in the derived category of graded right $A^{e}$-modules, for some graded algebra automorphism $\nu\in \Aut A$ called the Nakayama automorphism of $A$.  In the $\mathbb{Z}$-algebra case, it is unclear if AS-regularity is enough to guarantee the existence of such an isomorphism, and so we impose it as part of the definition of ASF$^{++}$-regularity.

The following is the main result of the paper, characterizing those noncommutative projective $\mathbb{Z}$-schemes associated to an ASF$^{++}$-regular $\ZZ$-algebra (Theorem \ref{thm.cs1}, Theorem \ref{thm.npl}, and Theorem \ref{thm.cs++2}), which extends \cite[Theorem 4.2]{BP}, \cite[Theorem 4.1]{MU1}, thus providing an answer to the $\mathbb{Z}$-algebra version of Question \ref{quest.triple}.  Before we state it, we remark that we will abuse notation in this paper by writing $\sC \cong \sD$ for categories $\sC$ and $\sD$ if they are {\it equivialent} (not necessarily isomorphic) categories.

\begin{theorem} \label{con.cs}
Let $\sC$ be a $k$-linear abelian category.
Then $\sC\cong \tails C$ for some right coherent ASF$^{++}$-regular $\ZZ$-algebra $C$ of dimension at least 1 and of Gorenstein parameter $\ell$ if and only if
\begin{enumerate}
\item[(GH1)] $\sC$ has a canonical bimodule $\omega _{\sC}$, and
\item[(GH2)] $\sC$ has an ample sequence $\{E_i\}_{i\in \ZZ}$ which is a full geometric helix  of period $\ell$ for $\sD^b(\sC)$.
\end{enumerate}
In fact, if $\sC$ satisfies (GH1) and (GH2), then $C:=C(\sC, \{E_i\}_{i\in \ZZ})_{\geq 0}$ is a right coherent ASF$^{++}$-regular $\ZZ$-algebra of dimension $\gldim \sC+1$ and of Gorenstein parameter $\ell$ such that $\sC\cong \tails C$.

Moreover, $C$ constructed above is right noetherian if and only if $\sC$ is a noetherian category.
\end{theorem}
The condition (GH1) requires that $\sC$ has an autoequivalence which induces a Serre functor on $\sD^{b}(\sC)$.  The notion of helix we use in (GH2) is similar to that of \cite{BP} (see \cite[Remark 3.17]{MU1}).  In comparison to \cite[Theorem 4.1]{MU1}, Theorem \ref{con.cs} is somewhat simpler in that no autoequivalence on $\sC$ (other than the canonical bimodule) is required.  Our proof of Theorem \ref{con.cs} requires the foundations of homological algebra for connected $\mathbb{Z}$-algebras developed in \cite{MN} and \cite{CMN}.  In particular, we use a variant of local duality \cite[Theorem 2.1]{CMN} in this paper to prove, in Theorem \ref{thm.serreduality}, that $\tails C$, where $C$ is ASF$^{++}$-regular, has a Serre functor.  Our argument is adapted from \cite[Appendix A]{NV}.

As an application of Theorem \ref{con.cs}, we construct a family of right noetherian AS-regular $\ZZ$-algebras from noncommutative quadric hypersurfaces.   In particular, we will show that every smooth quadric hypersurface and every standard noncommutative smooth quadric surface has a right noetherian AS-regular $\ZZ$-algebra as a homogeneous coordinate algebra (Theorem \ref{thm.nqh} and Theorem \ref{thm.a2}).

We now briefly describe the contents of the paper.  In Section 2 we recall relevant definitions and results from the theory of $\mathbb{Z}$-algebras we will need.  Although some of this material appears in \cite{MN}, most does not appear elsewhere, and is necessary for defining the notion of an ASF$^{++}$-regular $\ZZ$-algebra.  In Section 3, after recalling the notion of  noncommutative projective $\mathbb{Z}$-scheme, we prove variants of the version of local duality from \cite{CMN} which we will use in the proof of our main theorem.  We also include a number of results about derived functors and related triangles which will be used in the sequel.

In Section 4, we continue the study of regularity for $\mathbb{Z}$-algebras initiated in the paper \cite{MN}.  In \cite{MN}, we defined two notions of regularity for a $\ZZ$-algebra, namely, AS-regularity and ASF-regularity.  In this paper, after reviewing the aforementioned notions of regularity, we define two more notions of regularity for a $\ZZ$-algebra, namely ASF$^+$-regularity and ASF$^{++}$-regularity, and show the implications:
\begin{center}
\begin{tabular}{cccccccccccc}
ASF$^{++}$ & $\Rightarrow$ & ASF$^+$ & $\Rightarrow$ & AS & $\Rightarrow$ & ASF \\
\end{tabular}
\end{center}
 (Theorem \ref{thm.asasf}, Theorem \ref{thm.asfas2}, and Theorem \ref{thm.++}).  Recall that AS-regularity and ASF-regularity for a $\mathbb{Z}$-algebra $A$ are the same if $A$ has a ``balanced dualizing complex" \cite[Theorem 7.10]{MN}.  In this paper, we show that
ASF$^+$-regularity and ASF$^{++}$-regularity are the same if $A$ is $\ell$-periodic (Theorem \ref{thm.++}).  We note that the notion of ASF$^{++}$-regularity was first introduced and studied in \cite[Definition 3.1]{CMN}.

In Section 5, we introduce several concepts we will need for the proof of Theorem \ref{con.cs}, including various notions of helix.  The proof of Theorem \ref{con.cs}, as well as some consequences, are given in Section 6.  Finally, an application of our main result, to noncommutative quadric hypersurfaces, concludes the paper.

\section{$\ZZ$-algebras}

Although some of the results in this section are known (see \cite{MN}) or easy to see, we will discuss them rather carefully, hoping that this section (together with our previous paper \cite {MN}) will serve as a useful reference.

\subsection{$\mathbb{Z}$-algebras}

A {\it $\ZZ$-algebra} is an algebra with vector space decomposition $C=\oplus _{i, j\in \ZZ}C_{ij}$ and with the multiplication
$$C_{ij}\otimes C_{st}\to \begin{cases} C_{it} & \textnormal { if } j=s \\ 0 & \textnormal { if } j\neq s. \end{cases}$$  A $\ZZ$-algebra $C$ does not have a unity, but we assume that each subalgebra $C_{ii}$ has a unity $e_i\in C_{ii}$, called a local unity, so that $C_{ij}=e_iCe_j$ (and that $e_iae_j=a$ for every $a\in C_{ij}$).
Let $C, C'$ be $\ZZ$-algebras.  A {\it $\ZZ$-algebra homomorphism} $\phi:C\to C'$ is an algebra homomorphism $\phi:C\to C'$ such that $\phi (C_{ij})\subset C'_{ij}$ for all $i, j\in \ZZ$, and $\phi(e_i)=e_i'$ for all $i\in \ZZ$.
We say that $C$ is {\it locally finite} if $\dim _kC_{ij}<\infty$ for all $i, j$, and $C$ is {\it connected} if $C_{ij}=0$ for all $i>j$ and $C_{ii}=k$ for all $i$.

Let $C$ be a $\ZZ$-algebra.  A {\it graded right $C$-module} is a right $C$-module $M=\oplus _{j\in \ZZ}M_j$ with the action $M_i\otimes C_{ij}\to M_j$.  We assume that each $M_i$ is a unitary $C_{ii}$-module in the sense that $me_i=m$ for every $m\in M_i$.
The category of graded right $C$-modules is denoted by $\GrMod C$ whose morphisms are right $C$-module homomorphisms preserving degrees.  A {\it graded left $C$-module} is a left $C$-module $M=\oplus _{i\in \ZZ}M_i$ with the action $C_{ij}\otimes M_j\to M_i$.

Let $C, C'$ be $\ZZ$-algebras.  A {\it bigraded $C$-$C'$ bimodule} is a $C$-$C'$ bimodule $M=\oplus _{i,j}M_{ij}$ such that $e_iM:=\oplus _jM_{ij}$ is a graded right $C'$-module for every $i$ and $Me'_j:=\oplus _iM_{ij}$ is a graded left $C$-module for every $j$, that is, we have maps $M_{li}\otimes C'_{ij}\to M_{lj}$ and $C_{ij}\otimes M_{jl}\to M_{il}$.  A {\it homomorphism of bigraded $C$-bimodules} $M, N$ is a homomorphism $\phi:M\to N$ of $C$-bimodules such that $\phi(M_{ij})\subset N_{ij}$ for every $i, j\in \ZZ$.

For a graded left $C$-module $M$ and a graded right $C'$-module $N$, $M\otimes _kN:=\oplus _{i, j\in \ZZ}(M_i\otimes _kN_j)$ is naturally a bigraded $C$-$C'$ bimodule.
Note that $C$ itself is a bigraded $C$-bimodule.  If $C$ is connected, then $C_{\geq n}:=\oplus _{j-i\geq n}C_{ij}$ is a bigraded $C$-bimodule for every $n\in \NN$.

We define a graded right $C$-module
$$
P_i:=e_iC=\oplus _{j\in \ZZ}C_{ij}
$$
for every $i\in \ZZ$, and a graded left $C$-module
$$
Q_j:=Ce_j=\oplus _{i\in \ZZ}C_{ij}
$$
for every $j\in \ZZ$.  If $C$ is connected, then
$$
S_j:=e_j(C/C_{\geq 1})=(C/C_{\geq 1})e_j=P_je_j=e_jQ_j=e_jCe_j=C_{jj}
$$
has a structure of  a bigraded $C$-bimodule for every $j\in \ZZ$.  Moreover, $\{P_j\}_{j\in\ZZ}$ is the set of all indecomposable graded projective right $C$-modules up to isomorphism, and $\{S_j\}_{j\in \ZZ}$ is the set of all graded simple right $C$-modules up to isomorphism.  Note that $S_j$ is the unique simple quotient of $P_j$.

If $M$ is a graded right $C$-module, then $DM=\oplus _{i\in \ZZ}D(M_i)$ is a graded left $C$-module via $(af)(x)=f(xa)$ where $a\in C_{ij}, x\in M_j, f\in D(M_i)=(DM)_i$ so that $af\in D(M_j)=(DM)_j$.  Similarly, if $M$ is a graded left $C$-module, then $DM$ is a graded right $C$-module.  If $M$ is a bigraded $C$-$C'$ bimodule, then $DM:=\oplus _{i, j}D(M_{ji})$ is naturally a bigraded $C'$-$C$ bimodule.  It follows that
\begin{align*}
& D(e_iM)=\oplus _jD((e_iM)_j)=\oplus _jD(M_{ij})=\oplus _j(DM)_{ji}=(DM)e'_i \\
& D(Me'_j)=\oplus _iD((Me'_j)_i)=\oplus _iD(M_{ij})=\oplus _i(DM)_{ji}=e_j(DM)
\end{align*}
for every $i \in \ZZ$ and every $j \in \ZZ$, respectively.

We say that $M\in \GrMod C$ is {\it locally finite} if $\dim _kM_i<\infty$ for every $i$.  Note that $C$ is locally finite if and only if $P_j$ is locally finite for every $j$.  If $M\in \GrMod C$ is locally finite, then $DDM\cong M$ in $\GrMod C$.


The {\it opposite $\ZZ$-algebra of $C$} is the opposite algebra $C^o$ with $C^o_{ij}:=C_{-j, -i}$.  For $a^o\in C^o_{jk}=C_{-k, -j}, b^o\in C^o_{ij}=C_{-j, -i}$, $b^oa^o=ab \in C_{-k, -i}=C^o_{ik}$, so $C^o$ is in fact a $\ZZ$-algebra. If $M$ is a graded left $C$-module, then $\oplus_i M_{-i}$ is a graded right $C^o$-module under the action of $a^o\in C^o_{ij}=C_{-j, -i}$ on $x\in M_{-i}$ defined by $xa^o=ax\in M_{-j}$.

In fact, the category of graded left $C$-modules is equivalent to the category of graded right $C^o$-modules, so we often identify these two categories.  Similarly, we can see that the category of graded right $C$-modules is equivalent to the category of graded left $C^o$-modules.  Both of these facts are proven in \cite[Proposition 2.2]{MN}, where the algebra $C^o$ is denoted $\widetilde{C^{op}}$.

Note that if $C$ is a connected $\ZZ$-algebra, then $C^o$ is again a connected $\ZZ$-algebra.  Let $C, C'$ be $\ZZ$-algebras. The category of bigraded $C$-$C'$ bimodules is denoted by $\Bimod (C - C')$.  Note that the categories ${\GrMod C}$ and $\Bimod (C - C')$ are Grothendieck categories \cite[Section 3]{Vq}, \cite[Proposition 2.2(1)]{MN}, hence abelian.

\begin{remark} \label{rem.opo}  Let $C$ be a $\ZZ$-algebra.
\begin{enumerate}
\item{} If $M=\oplus _{i\in \ZZ}M_i$ is a graded left $C$-module, then, precisely speaking, $\oplus _{i\in \ZZ}M_{-i}\in \GrMod C^o$ is a graded right $C^o$-module, however, we often identify them in this paper.
\item{} If $M$ is a bigraded $C$-bimodule, then $e_iM:=\oplus _{j\in \ZZ}M_{ij}$ is a graded right $C$-module for every $i\in \ZZ$, and $Me_j:=\oplus _{i\in \ZZ}M_{ij}$ is a graded left $C$-module for every $j\in \ZZ$, but $Me_j$
is not a graded right $C^o$-module in this grading.  By defining $M^o:=\oplus _{ij}M_{-j, -i}\in \Bimod (C^o-C^o)$, we identify $Me_j$ with a graded right $C^o$-module $e_{-j}^oM^o:=\oplus _{i\in \ZZ}M^o_{-j, i}:=\oplus _{i\in \ZZ}M_{-i, j}=\oplus _{i\in \ZZ}(Me_{j})_{-i}$.

\item{} We write
\begin{align*}
P_i^o & :=e_i^oC^o=\oplus _{j\in \ZZ}C^o_{ij}=\oplus _{j\in \ZZ}C_{-j, -i}=Ce_{-i}=:Q_{-i}, \\
Q_j^o & :=C^oe_j^o=\oplus _{i\in \ZZ}C^o_{ij}=\oplus _{i\in \ZZ}C_{-j, -i}=e_{-j}C=:P_{-j}, \\
S_j^o & =e_j^oC^oe_j^o=C^o_{jj}=C_{-j, -j}=e_{-j}Ce_{-j}=:S_{-j}.
\end{align*}
In particular, we identify a graded right $C^o$-module $P_j^o:=e^o_jC^o$ with a graded left $C$-module $Q_{-j}:=Ce_{-i}$, and so on.
\item{} $D:\GrMod C\to \GrMod C^o$ is defined by $D(\oplus _iM_i):=\oplus_i D(M_{i})$ when we view $D(\oplus _iM_i)$ as a graded left $C$-module while $D(\oplus _iM_i):=\oplus_i D(M_{-i})$ when we view $D(\oplus _iM_i)$ as a graded right $C^o$-module.

\end{enumerate}

\end{remark}

\subsection{$\Hom$ and $\otimes$}

If $M$ is a bigraded $C'$-$C$ bimodule and $N$ is a graded right $C$-module, then $\uHom_C(M, N):=\oplus _i\Hom_C(e_iM, N)=\oplus _i\Hom_C(\oplus _kM_{ik}, N)$ has a structure of a graded right $C'$-module
via
$(fa)(x)=f(ax)$ where $a\in C'_{\ell i}, x\in \oplus _jM_{ij}, f\in \uHom_C(M, N)_\ell=\Hom_C(\oplus _jM_{lj}, N)$ by the map $a\cdot :M_{ij}\to M_{\ell j}$  so that $fa\in \uHom_C(M, N)_\ell=\Hom_C(\oplus _jM_{\ell j}, N)$.  (Although $\oplus _jM_{ij}$ does not have a structure of a left $C'$-module, the left multiplication $a\cdot :\oplus _jM_{ij}\to \oplus _j M_{\ell j}$ is well-defined, which induces the right action $\cdot a:\Hom_C(\oplus _jM_{ij}, N)\to \Hom_C(\oplus _jM_{\ell j}, N)$.)

If $M$ is a graded right $C$-module and $N$ is a bigraded $C'$-$C$ bimodule, then $\uHom_C(M, N):=\oplus _j\Hom_C(M, e_jN)=\oplus _j\Hom_C(M, \oplus _\ell N_{j\ell})$ has a structure of a graded left $C'$-module via $(af)(x)=a(f(x))$ where $a\in C'_{ij}, x\in M, f\in \uHom_C(M, N)_j=\Hom_C(M, \oplus _\ell N_{j\ell})$ by the map $a\cdot :N_{j\ell}\to N_{i\ell}$ so that $af\in \uHom_C(M, N)_i=\Hom_C(M, \oplus _\ell N_{i\ell})$.

If $M$ is a bigraded $C'$-$C$ bimodule and $N$ is a bigraded $C''$-$C$ bimodule,
then $\uHom_C(M, N):=\oplus _{i,j}\Hom_C(e_jM, e_iN)$ has a structure of a bigraded $C''$-$C'$ bimodule since
\begin{eqnarray*}
e''_i\uHom_C(M, N) & = & \oplus _{j\in \ZZ}\uHom_C(M, N)_{ij} \\
& = & \oplus _{j\in \ZZ}\Hom_C(e'_jM, e''_iN) \\
& = & \uHom_C(M, e''_iN)
\end{eqnarray*}
is a graded right $C'$-module for every $i\in \ZZ$, and
\begin{eqnarray*}
\uHom_C(M, N)e'_j & = & \oplus _{i\in \ZZ}\uHom_C(M, N)_{ij} \\
& = & \oplus _{i\in \ZZ}\Hom_C(e'_jM, e''_iN) \\
& = & \uHom_C(e'_jM, N)
\end{eqnarray*}
is a graded left $C''$-module for every $i\in \ZZ$.

The proof of the following lemma is straightforward and omitted.

\begin{lemma} \label{lem.Hom}
Let $C$ be a $\ZZ$-algebra.
\begin{enumerate}
\item{} For $M\in \GrMod C$, $\uHom_C(C, M)\cong M$ in $\GrMod C$ so that
$$
\Hom_C(P_i, M)=M_i
$$
for every $i\in \ZZ$. In particular,
$$\Hom_C(P_i, S_j)=\begin{cases} S_j & \textnormal { if } i=j \\
0 & \textnormal { if } i\neq j.\end{cases}$$
\item{} $\uHom_C(P_i, C)\cong Q_i$ in $\GrMod C^o$,
and $\uHom_{C^o}(Q_i, C^o)\cong P_i$ in $\GrMod C$ for every $i\in \ZZ$.
\end{enumerate}
\end{lemma}

\begin{definition} Let $C, C', C''$ be $\ZZ$-algebras.  For $M\in \GrMod C$ and $N\in \GrMod C^o$, we define
$$M\otimes _CN:=\Coker(\oplus _{i, j\in \ZZ}(M_i\otimes _{C_{ii}}C_{ij}\otimes _{C_{jj}}N_j)\to \oplus _{k\in \ZZ}M_k\otimes _{C_{kk}}N_k)$$
where the morphism is induced by the usual difference between left and right multiplication.

\begin{itemize}
\item{} For $M\in \GrMod C$ and $N\in \Bimod (C - C')$, we define
$$
M\uotimes _CN=\oplus _{\ell \in \ZZ}(M\otimes _CNe_{\ell}')\in \GrMod C.
$$
\item{} For $M\in \Bimod (C' - C)$ and $N\in \GrMod C^o$, we define
$$
M\uotimes _CN=\oplus _{i\in \ZZ}(e_iM\otimes _CN)\in \GrMod {C'}^o.
$$
\item{} For $M\in \Bimod (C - C'), N\in \Bimod (C' - C'')$, we define
$$
M\uotimes _CN=\oplus _{i, j\in \ZZ}(e_iM\otimes _CNe''_j)\in \Bimod (C - C'').
$$
\end{itemize}
\end{definition}
Note that from general properties of adjoint functors and \cite[Proposition 5.3(2)]{MN}, for $N\in \Bimod (C' - C'')$, the functor
$$
\uotimes_C N:\Bimod (C - C') \rightarrow \Bimod (C - C'')
$$
commutes with colimits.

\begin{lemma}
\cite[Section 4.1, Proposition 5.3]{MN}
Let $C$ be a $\ZZ$-algebra.
\begin{enumerate}
\item{} For $M\in \GrMod C$, $M\uotimes _CC\cong M$ in $\GrMod C$ so that $M\uotimes _CQ_j\cong M_j$ for every $j\in \ZZ$.
\item{} For $N\in \GrMod C^o$, $C\uotimes _CN\cong N$ in $\GrMod C^o$ so that $P_i\uotimes _CN\cong N_i$ for every $i\in \ZZ$.
\item{} For $M\in \GrMod C, N\in \GrMod C'$ and $L\in \Bimod (C - C')$,
$$\Hom_{C'}(M\uotimes_CL, N)\cong \Hom_C(M, \uHom_{C'}(L, N)).$$
\end{enumerate}
\end{lemma}

\subsection{Noetherian and Coherent Properties}

For a set of objects $\cE$ in an additive category $\sC$, we denote by $\add \cE$ the set of objects in $\sC$ consisting of all finite direct sums of objects in $\cE$.

\begin{remark}

The notation $\add \cE$ usually denotes the set of objects in
$\sC$ consisting of direct summands of finite direct sums of objects in $\cE$, so the above notation is not standard.

\end{remark}

\begin{definition} Let $C$ be a $\ZZ$-algebra.
\begin{enumerate}
\item{} We say that $M\in \GrMod C$ is {\it finitely generated} (resp. {\it finitely presented}) if there exists an exact sequence $F^0\to M\to 0$ (resp. $F^1\to F^0\to M\to 0$) where $F^i \in \add \{P_{j}\}_{j \in \mathbb{Z}}$.

\item{} We say that $M\in \GrMod C$ is {\it coherent} if $M$ is finitely generated and $\Ker \phi$ is finitely generated for every homomorphism $\phi: F\to M$ with $F \in \add \{P_{j}\}_{j \in \mathbb{Z}}$.
\item{} We denote by $\grmod C$ the full subcategory of $\GrMod C$ consisting of finitely presented modules, and by $\coh C$ the full subcategory of $\GrMod C$ consisting of coherent modules.
\item{} We say that $C$ is {\it right coherent} if $P_i, S_i\in \coh C$ for every  $i\in \ZZ$.
\end{enumerate}
\end{definition}

We call a module in $\add \{P_{j}\}_{j \in \mathbb{Z}}$ {\it finitely generated free}.
In this terminology, $C=\oplus _jP_j\in \GrMod C$ itself is free but not finitely generated free.  If $C$ is connected, then every finitely generated projective graded right $C$-module is isomorphic to a module in $\add \{P_{j}\}_{j \in \mathbb{Z}}$.

The following result will often be used without comment in the sequel.
\begin{lemma} \label{lem.ababc}
Let $C$ be a locally finite connected $\ZZ$-algebra.
\begin{enumerate}
\item{} $\coh C$ is an abelian category.
\item{} If $C$ is right coherent, then $\grmod C=\coh C$ so that $\grmod C$ is an abelian category.
\item{} Conversely, if $\grmod C$ is an abelian category, then $P_j\in \coh C$ for every $j\in \ZZ$.
\end{enumerate}
\end{lemma}

\begin{proof} (1) This follows from \cite[Proposition 1.1]{P}.

(2) If $M\in\coh C$, then $M\in \grmod C$ by definition.   Conversely, if $C$ is right coherent and $M\in \grmod C$, then there exists an exact sequence $F^1\to F^0\to M\to 0$ in $\GrMod C$ where $F^1, F^0\in \add \{P_{j}\}_{j \in \mathbb{Z}} \subset \coh C$.  Since $\coh C$ is an abelian category by (1), $M\cong \Coker(F^1\to F^0)\in \coh C$.

(3) If $\grmod C$ is an abelian category, then, for every homomorphism $\phi:F\to P_j$ with $F, P_j\in \add \{P_{j}\}_{j \in \mathbb{Z}} \subset \grmod C$,
$\Ker \phi\in \grmod C$.  In particular, $P_j, \Ker \phi$ are finitely generated, so $P_j\in \coh C$.
\end{proof}

\begin{definition} We say that a $\ZZ$-algebra $C$ is right noetherian if $P_j\in \GrMod C$ is a noetherian object for every $j\in \ZZ$.
\end{definition}

By \cite[Definition 3.1]{Vq}, $C$ is right noetherian if and only if $\GrMod C$ is a locally noetherian (Grothendieck) catetgory.

\begin{lemma}
\label{lem.cone}
If $C$ is a right noetherian $\ZZ$-algebra, then every noetherian module is coherent.  In particular, $C$ is right coherent and $\grmod C$ is an abelian category.
\end{lemma}

\begin{proof} If $M\in \GrMod C$ is a noetherian module, then $M$ is finitely generated, so there exists $F\in \add \{P_{j}\}_{j \in \mathbb{Z}}$ and a surjection $\phi:F\to M$.  Since $P_j$ is noetherian for every $j\in \ZZ$, $F$ is noetherian, so $\Ker \phi\subset F$ is noetherian.  It follows that $\Ker \phi$ is finitely generated, so $M\in \coh C$.  In particular, since $P_j, S_j\in \GrMod C$ are noetherian for every $j\in \ZZ$, $P_j, S_j\in \coh C$, so $C$ is right coherent.
\end{proof}

\begin{lemma} \label{lem.noe}
A $\ZZ$-algebra $C$ is right noetherian if and only if $\grmod C$ is a noetherian category.
\end{lemma}

\begin{proof}
For every $M\in \grmod C$, there exists a surjection $F\to M$ in $\GrMod C$ where $F\in \add \{P_{j}\}_{j \in \mathbb{Z}}\subset \GrMod C$ is a noetherian object, so $M\in \GrMod C$ is a noetherian object, hence $\grmod C$ is a noetherian category.

Conversely, if $\grmod C$ is a noetherian category, then $P_j\in \grmod C$ is a noetherian object for every $j\in \ZZ$, so $C$ is right noetherian.
\end{proof}

\subsection{Module Categories}

Let $C$ be a $\ZZ$-algebra.  We denote by $\Mod C$ the category of right $C$-modules which is ``unitary" in the sense that $M=MC$.
The following lemma is known (cf. \cite{S}, \cite{Vq}).  We give a proof for the convenience of the reader.

\begin{lemma} \label{lem.mgm} For every $\ZZ$-algebra $C$, $\GrMod C\cong \Mod C$.
\end{lemma}

\begin{proof}
Let $M\in \Mod C$.  Since $e_ie_j=\begin{cases} e_j & \textnormal { if } i=j \\ 0 & \textnormal { if } i\neq j, \end{cases}$ if $me_i=ne_j\in Me_i\cap Me_j$ for $i\neq j$, then $ne_j=ne_j^2=me_ie_j=0$, so $\oplus _{j\in \ZZ}Me_j\subset M$.  Since $M=MC$, for every $m\in M$, there exist $n\in M$ and $a=\sum_{i, j}a_{ij}\in \oplus _{i, j\in \ZZ}C_{ij}=C$ such that $m=na=\sum _{i, j}na_{ij}=\sum _j(\sum_ina_{ij})e_j\in \oplus _{j\in \ZZ}Me_j$, so $M=\oplus _{j\in \ZZ}Me_j$.  It is easy to see that $\oplus _{i\in \ZZ}Me_i$ is naturally a graded right $C$-module.  If $\phi:M\to N$ is a homomorphism of right $C$-modules, then $\phi (me_j)=\phi (m)e_j$, so $\phi (Me_j)\subset Ne_j$, hence $\phi$ is naturally a homomorphism of graded right $C$-modules.  It follows that
 $\Mod C\to \GrMod C; \; M\mapsto \oplus _{j\in \ZZ}Me_j$ is an equivalence functor.
\end{proof}

\begin{remark}  By Lemma \ref{lem.mgm}, we have the following:
\begin{enumerate}
\item{} For a $\ZZ$-algebra $C$ and $M, N\in \GrMod C$, $M\cong N$ in $\GrMod C$ if and only if $M\cong N$ in $\Mod C$.
\item{} For $\ZZ$-algebras $C, C'$,  $\GrMod C\cong \GrMod C'$ if and only if $\Mod C\cong \Mod C'$.
\end{enumerate}
\end{remark}

\begin{lemma} \label{lem.mgm2} Let $C, C'$ be $\ZZ$-algebras.  If $C\cong C'$ as algebras (not necessarily as $\ZZ$-algebras), then $\GrMod C\cong \GrMod C'$.
\end{lemma}

\begin{proof}
 Let $\phi:C\to C'$ be an isomorphism of algebras.  If $M\in \Mod C'$, then $MC'=M$, so $M_{\phi}C:=M\phi(C)=MC'=M$, hence $M_{\phi}\in \Mod C$.  It follows that $\GrMod C\cong \Mod C\cong \Mod C'\cong \GrMod C'$ by Lemma \ref{lem.mgm}.
\end{proof}

\subsection{Periodicity}

Let $C, C'$ be $\ZZ$-algebras, and $\phi:C\to C'$ a homomorphism of $\ZZ$-algebras.  For $M\in \GrMod C'$,
we define $M_{\phi}\in \GrMod C$ by $M_{\phi}=M$ as a graded vector space with the action $m*a=m\phi(a)$ for $m\in M, a\in C$.  This induces a functor $(-)_{\phi}:\GrMod C'\to \GrMod C$.

For a bigraded vector space $M$ and $r, \ell\in \ZZ$, we define a bigraded vector space $M(r, \ell)$ by $M(r, \ell)_{ij}=M_{r+i, \ell+j}$.

\begin{lemma} \label{lem.rel}
Let $C, C'$ be $\ZZ$-algebras.
\begin{enumerate}
\item{} For $\ell\in \ZZ$, $C(\ell):=C(\ell, \ell)$ is a $\ZZ$-algebra.
\item{} If $M$ is a bigraded $C$-$C'$ bimodule, then $M(\ell, r)$ is a bigraded $C(\ell)$-$C'(r)$ bimodule.
\item{} If $\phi:C\to C(-\ell)$ is an isomorphism of $\ZZ$-algebras, then $\phi^{-1}$ induces an isomorphism $C\to C(\ell)$ of $\ZZ$-algebras, which is denoted by $\phi^{-1}$ by abuse of notation, and $C(0,-\ell)_{\phi}\cong {_{\phi^{-1}}C(\ell, 0)}$
as bigraded $C$-bimodules.
\end{enumerate}
\end{lemma}

\begin{proof} (1) Note that $C=C(\ell)$ as ungraded algebras.  Since
$$C(\ell)_{ij}\otimes C(\ell)_{jk}=C_{i+\ell, j+\ell}\otimes C_{j+\ell, k+\ell}\to C_{i+\ell, k+\ell}=C(\ell)_{ik},$$
the result follows.

(2) Note that $M(\ell, r)=M$ as ungraded $C$-$D$ bimodules.  Since
\begin{align*}
& C(\ell)_{ij}\otimes M(\ell, r)_{jk}=C_{\ell+i, \ell+j}\otimes M_{\ell+j, r+k}\to M_{\ell+i, r+k}=M(\ell, r)_{ik}, \\
& M(\ell, r)_{ij}\otimes D(r)_{jk}=M_{\ell+i, r+j}\otimes D_{r+j. r+k}\to M_{\ell+i, r+k}=M(\ell, r)_{ik},
\end{align*}
the result follows.

(3) Since $\phi^{-1}:C\to C(\ell)$ is an isomorphism of (ungraded) algebras, and
$$\phi^{-1}(C_{ij})=\phi^{-1}(C(-\ell)_{i+\ell, j+\ell})=C_{i+\ell, j+\ell}=C(\ell)_{ij},$$
$\phi^{-1} :C\to C(\ell)$ is an isomorphism of $\ZZ$-algebras.

Since $\phi:{_{\phi^{-1}}C(\ell,0)}\to C(0,-\ell)_{\phi}$ is an isomorphism of bigraded vector spaces, and

$$\phi (a*mb)=\phi(\phi^{-1}(a)mb)=a\phi(m)\phi(b)=a\phi(m)*b,$$
for $m\in {_{\phi^{-1}}C(0, -\ell)}$ and $a, b\in C$, $\phi:{_{\phi^{-1}}C(\ell,0)}\to C(0,-\ell)_{\phi}$ is an isomorphism of bigraded $C$-bimodules.
\end{proof}

Let $C$ be a $\ZZ$-algebra.  Since $C\cong C(r)$ as algebras (but not necessarily as $\ZZ$-algebras in general), $\GrMod C\cong \GrMod C(r)$ by Lemma \ref{lem.mgm2}.  In fact, we have an equivalence functor $(r):\GrMod C\to \GrMod C(r)$ defined by $M(r):=\oplus _{j\in \ZZ}M_{r+j}$.  Since $(e_iC)(r)=\oplus _{j\in \ZZ}C_{i, j+r}=e_{i-r}(C(r))$, the assignment $F\mapsto F(r)$ preserves finitely generated projectives, so $(r)$ induces an equivalence functor $(r):\grmod C\to \grmod C(r)$.
For $M\in \GrMod C$, $M(r)$ is not a graded right $C$-module in general, so there exists no notion of homomorphism of graded right $C$-modules of degree $r$.

\begin{lemma} \label{lem.rell}
Let $C, C', C''$ be $\ZZ$-algebras.
\begin{enumerate}
\item{} For $\ell\in \ZZ$, $-\uotimes_CC(0, \ell)\cong (-)(\ell):\GrMod C\to \GrMod C(\ell)$ as functors.
\item{} If $\phi:C''\to C'$ is a homomorphism of $\ZZ$-algebras and $M$ is a bigraded $C$-$C'$ bimodule, then $-\uotimes_CM_{\phi}\cong (-\uotimes _CM)_{\phi}:\GrMod C\to \GrMod C''$ as functors.
\item{} In particular, if $\phi:C\to C(\ell)$ is a homomorphism of $\ZZ$-algebras, then $-\otimes_CC(0, \ell)_{\phi}\cong (-)(\ell)_{\phi}:\GrMod C\to \GrMod C$ as functors.
\end{enumerate}
\end{lemma}

\begin{proof} Functors of both sides are naturally isomorphic on ungraded module categories and they are compatible with the grading, hence the result.
\end{proof}

We say that $C$ is {\it $r$-periodic} if $C(r)\cong C$ as $\ZZ$-algebras.   If $C$ is $r$-periodic and $\phi:C\to C(r)$ is an isomorphism of $\ZZ$-algebras, then, for $M\in \GrMod C$, $M(r)_{\phi}\in \GrMod C$, so there exist autoequivalences of $\GrMod C$ and $\grmod C$ defined by $M\mapsto M(r)_{\phi}$.

\begin{lemma} \label{lem.pper}
If $C$ is an $r$-periodic $\ZZ$-algebra, then an isomorphism $\phi:C\to C(r)$ of $\ZZ$-algebras restricts to an isomorphism $P_\ell \to P_{\ell+r}(r)_{\phi}$ in $\GrMod C$ for every $\ell \in \ZZ$.
\end{lemma}

\begin{proof} If $\phi:C\to C(r)$ is an isomorphism of $\ZZ$-algebras, then $P_{\ell+r}(r)_{\phi}=\oplus _iC_{\ell+r, i+r}=\phi (P_\ell)$ as graded vector spaces.  For $a\in P_{\ell+r}(r)_{\phi}$ and $b\in C$,  $\phi(a)*b=\phi(a)\phi(b)=\phi(ab)$, so we have a commutative diagram
$$\begin{CD}
(P_\ell)_i\otimes C_{ij}= & C_{\ell i}\otimes C_{ij} @>\cdot>> C_{\ell j} & =(P_\ell)_j \\
& @V{\phi\otimes \id}VV @VV\phi V \\
P_{\ell+r}(r)_i\otimes C_{ij}= & C_{\ell+r, i+r}\otimes C_{ij} @>*>> C_{\ell+r, j+r} & =P_{\ell+r}(r)_j
\end{CD}$$
hence $\phi:P_\ell\to P_{\ell+r}(r)_{\phi}$ is an isomorphism in $\GrMod C$ for every $\ell \in \ZZ$.
\end{proof}

\subsection{$\ZZ$-algebras Associated to Graded Algebras}

For a graded algebra $A$, we define a $\ZZ$-algebra $\overline A$ by $\overline A:=\oplus _{i, j\in \ZZ}A_{j-i}$.
The following lemma is well-known (cf. \cite{BP}, \cite{P}, \cite{S}).

\begin{lemma} \label{lem.ovwi}
For a graded algebra $A$, the functors $\Phi:\GrMod A\to \GrMod \overline A$ and $\Phi:\grmod A\to \grmod \overline A$ defined by $\Phi(M):=\oplus _{i\in \ZZ}M_{i}$ are equivalences of categories sending $A(-i)$ to $P_i$ and $k(-i)$ to $S_i$ for every $i\in \ZZ$.
\end{lemma}

\begin{proof} It is well-known that the functor $\Phi:\GrMod A\to \GrMod \overline A$ defined as above is an equivalence functor (cf. \cite{S}).  Since
\begin{align*}
\Phi(A(-i)) & :=\oplus _jA(-i)_{j}=\oplus _jA_{j-i}=\oplus _j\overline A_{ij}=e_i\overline A=:P_i, \\
\Phi(k(-i)) & :=\oplus _jk(-i)_{j}=\oplus _jk_{j-i}=e_i\overline Ae_i=:S_i,
\end{align*}
$\Phi$ restricts to an equivalence functor $\Phi:\grmod A\to \grmod \overline A$.
\end{proof}

\begin{lemma} \label{lem.asfz1} Let $A$ be a graded algebra.  Viewing
$\GrMod \overline A^o$ as the category of graded left $\overline A$-modules, the following hold.
\begin{enumerate}
\item{}
$\Phi^o:\GrMod A^o\to \GrMod \overline A^o$ defined by $\Phi^o(M):=\oplus _{j\in \ZZ}M_{-j}$ is an equivalence functor such that $\Phi^o(A(i))\cong \overline Ae_{i}$ for every $i\in \ZZ$.
\item{} $\GrMod A^e\to \Bimod (\overline A - \overline A); \; M\mapsto \overline M:=\oplus_{i, j\in \ZZ}M_{j-i}$ is a functor such that, for every $M\in \GrMod A^e$,  $\Phi(M(-i))\cong e_i\overline M$ for every $i\in \ZZ$, and $\Phi^o(M(i))\cong \overline Me_{-i}$ for every $i\in \ZZ$.   (Here we view $\Phi:\GrMod A^e\to \GrMod A\to \GrMod \overline A$ and $\Phi^o:\GrMod A^e\to \GrMod A^o\to \GrMod \overline A^o$ by composing with the natural functors $\GrMod A^e\to \GrMod A$ forgetting the graded left $A$-module structure and $\GrMod A^e\to \GrMod A^o$ forgetting the graded right $A$-module structure.)

\item{} $\Phi^oD=D\Phi:\GrMod A\to \GrMod \overline A^o$ as functors, and $\Phi D= D\Phi^o:\GrMod A^o\to \GrMod \overline A$ as functors.
\end{enumerate}
\end{lemma}

\begin{proof} (1) By Lemma \ref{lem.ovwi},
$\Phi^o:\GrMod A^o\to \GrMod \overline A^o$ defined by $\Phi^o(M):=\oplus _{j\in \ZZ}M_{j}$ is an equivalence functor such that $\Phi^o(A^o(-i))\cong e^o_{i}\overline A^o$ for every $i\in \ZZ$ viewed as a graded right $\overline A^o$-module, however, $\Phi^o(M)=\oplus _{j\in \ZZ}M_{-j}$ and $e^o_{i}\overline A^o=\overline Ae_{-i}$ viewed as a graded left $\overline A$-module by Remark \ref{rem.opo}.

(2) Let $M\in \GrMod A^e$.  For $x\in \overline M_{ij}=M_{j-i}, a\in \overline A_{si}=A_{i-s}, b\in \overline A_{jt}=A_{t-j}$, $axb\in M_{t-s}=\overline M_{st}$, so we may view $\overline M\in \Bimod (\overline A - \overline A)$.

Let $\phi\in \Hom_{A^e}(M, N)$ where $M, N\in \GrMod A^e$.  For $x\in \overline M_{ij}=M_{j-i}$, $\phi(x)\in N_{j-i}=\overline N_{ij}$, so we may define a map $\overline \phi:\overline M\to \overline N$ by $\overline \phi(x):=\phi(x)$ such that $\phi(\overline M_{ij})\subset \overline N_{ij}$.   Since $\overline \phi(axb)=\phi(axb)=a\phi(x)b=a\overline \phi(x)b$, $\overline \phi\in \Hom_{\overline A}(\overline M, \overline N)$.  Moreover, we have
$$\Phi(M(-i)):=\oplus _{j\in \ZZ}M(-i)_{j}=\oplus _{j\in \ZZ}M_{j-i}=\oplus _{j\in \ZZ}\overline M_{ij}=:e_i\overline M,$$
and
$$\Phi^o(M(i)):=\oplus _{j\in \ZZ}M(i)_{-j}=\oplus _{j\in \ZZ}M_{i-j}=\oplus _{j\in \ZZ}\overline M_{-j,-i}=:\overline Me_{-i}.$$

(3) Recall that $D:\GrMod A\to \GrMod A^o$ is defined by $D(\oplus _iM_i):=\oplus_i D(M_{-i})$ while $D:\GrMod \overline A\to \GrMod \overline A^o$ is defined by $D(\oplus _iM_i):=\oplus_i D(M_{-i})$ when we view $D(\oplus _iM_i)$ as a graded left $\overline A$-module (see Remark \ref{rem.opo} (4)).  Since
$$\Phi^oD(\oplus _iM_i)=\Phi^o(\oplus_iD(M_{-i}))=\oplus _iD(M_{i})=D(\oplus _iM_{i})=D\Phi(\oplus _iM_i)$$
in $\GrMod \overline A^o$ for $M\in \GrMod A$,
we see that $\Phi^oD= D\Phi:\GrMod A\to \GrMod \overline {A}^o$.

Similarly, since
$$\Phi D(\oplus _iM_i)=\Phi(\oplus_iD(M_{-i}))=\oplus _iDM_{-i}=D(\oplus _iM_{-i})=D\Phi^o(\oplus _iM_i)$$
in $\GrMod \overline A$ for $M\in \GrMod A^o$,
we see that $\Phi D= D\Phi^o:\GrMod A^o\to \GrMod \overline {A}$.
\end{proof}

\begin{remark} \label{rem.badu}
For a graded algebra $A$, $\overline {(A^o)}_{ij}={A^o}_{j-i}=A_{-i+j}=\overline A_{-j, -i}={(\overline A)^o}_{ij}$ for every $i, j\in \ZZ$, so $\overline {(A^o)}=(\overline A)^o$ as $\ZZ$-algebras.   In particular, $A$ is a connected graded algebra if and only if $A^o$ is a connected graded algebra if and only if $\overline A$ is a connected $\ZZ$-algebra if and only if $(\overline A)^o$ is a connected $\ZZ$-algebra.
\end{remark}

\begin{lemma} \label{lem.1p}  \cite [Proposition 3.1]{S} Let $C$ be a $\ZZ$-algebra.  Then there exists a graded algebra $A$ such that $C\cong \overline A$ as $\ZZ$-algebras if and only if $C$ is 1-periodic.
\end{lemma}

\section{Derived Categories of Graded Modules} \label{section.ld}

One of the main results in \cite{MN} and \cite{CMN} is a version of local duality for connected $\ZZ$-algebras.  In this section, we prove another version of local duality for $\ZZ$-algebras (Theorem \ref{thm.ld0}) and we prove other results about derived functors associated to noncommutative projective $\ZZ$-schemes which will be employed in the sequel.

Throughout the remainder of the paper, if $\sC$ denotes an abelian category, we let $\sD(\sC)$ (resp. $\sD^{-}(\sC)$, $\sD^{+}(\sC)$, $\sD^b(\sC)$) denote the derived category of $\sC$ (resp. the bounded above derived category of $\sC$, the bounded below derived category, the bounded derived category).  We will also utilize various left and right derived functors of functors already introduced, and the reader may consult \cite[Section 6]{MN} for more information about these derived functors.

\subsection{Noncommutative Projective $\ZZ$-schemes}

\begin{definition} Let $C$ be a $\ZZ$-algebra and $X\in \sD^b(\GrMod C)$.  A complex $(F, d)$
$$\begin{CD}\cdots @>d^2>> F^2 @>d^1>> F^1 @>d^0>> F^0 @>>> M @>>> 0\end{CD}$$
where $F^q\in \add \{P_{j}\}_{j \in \mathbb{Z}}$ is called a {\it finitely generated free resolution of $X$} if $F$ is quasi-isomorphic to $X$.  A finitely generated free resolution of $X$ is called {\it minimal} if $\Im d^q\subset F^qC_{\geq 1}$ for every $q\in \ZZ$.
\end{definition}

\begin{lemma} \label{lem.ab}
If $C$ is a right coherent connected $\ZZ$-algebra, then every $M\in \grmod C$ has a unique minimal finitely generated free resolution up to isomorphism.

\end{lemma}

\begin{proof}
The existence of a minimal finitely generated free resolution follows from \cite[Proposition 4.4]{MN}.
Using the same argument as in the connected graded case,
a minimal finitely generated free resolution over a connected $\ZZ$-algebra is unique up to isomorphism.
\end{proof}

The following condition is essential for the rest of the paper.

\begin{definition} A connected $\ZZ$-algebra $C$ is called {\it right $\Ext$-finite} if every $S_i$ has a minimal  finitely generated free resolution in $\GrMod C$.
\end{definition}

\begin{lemma} \label{lem.colf}
If $C$ is a right coherent connected $\ZZ$-algebra, then $C$ is right $\Ext$-finite.

\end{lemma}

\begin{proof} This follows from Lemma \ref{lem.ab}.
\end{proof}

\begin{lemma} \label{lem.lofi}
Let $C$ be a connected $\ZZ$-algebra.
If $C$ is right $\Ext$-finite, then $C$ is locally finite.  In particular, if $C$ is right coherent, then $C$ is locally finite.
\end{lemma}

\begin{proof} This follows from \cite[Remark 3.3]{MN}.
\end{proof}

In summary, we have the following implications for a connected $\ZZ$-algebra.
\begin{center} right noetherian $\Rightarrow$ right coherent $\Rightarrow$ right $\Ext$-finite $\Rightarrow$ locally finite. \end{center}

\begin{definition}
Let $C$ be a $\ZZ$-algebra.  We say that $M\in \GrMod C$ is {\it right bounded} if $M_{\geq m}=0$ for some $m\in \ZZ$ and {\it left bounded} if $M_{\leq m}=0$ for some $m\in \ZZ$.

We denote by $\Tors C$ the full subcategory of $\GrMod C$ consisting of modules $M$ such that $xC$ is right bounded for every $x\in M$.

We define a torsion functor $\t:\GrMod C\to \Tors C\subset \GrMod C$ by
$$\t(M):=\{x\in M\mid xC \textnormal { is right bounded}\}.$$
\end{definition}

\begin{lemma} \label{lem.loc}
If $C$ is a right $\Ext$-finite connected $\ZZ$-algebra, then
$\Tors C$ is a localizing subcategory of $\GrMod C$.
\end{lemma}

\begin{proof} By \cite [Lemma 3.5]{MN}, $\Tors C$ is a Serre subcategory of $\GrMod C$.  Since
$\t(M)$ is the largest torsion submodule of $M\in \GrMod C$,
$\Tors C$ is a localizing subcategory of $\GrMod C$ by \cite[Proposition 4.5.2]{Po}.
\end{proof}


\begin{definition} \label{def.nps}
Let $C$ be a right $\Ext$-finite connected $\ZZ$-algebra.  We define the quotient category  $\Tails C:=\GrMod C/\Tors C$.  We call $\Tails C$ the {\it noncommutative projective $\ZZ$-scheme associated to $C$}.
\end{definition}

Let $C$ be a right $\Ext$-finite connected $\ZZ$-algebra.  We denote by
$$
\pi : \GrMod C\to \Tails C
$$
the quotient functor.   Since $\Tors C$ is a localizing subcategory of $\GrMod C$ by Lemma \ref{lem.loc},  $\pi :\GrMod C\to \Tails C$ has a right adjoint
$$
\omega :\Tails C\to \GrMod C.
$$
We write
$$
Q:=\omega\pi:\GrMod C\to \GrMod C.
$$
We often write $\cM:=\pi M\in \Tails C$ for $M\in \GrMod C$.  We also write $\Ext^q_{\cC}(\cM, \cN):=\Ext^q_{\Tails C}(\cM, \cN)$ for $\cM, \cN\in \Tails C$.

If $C$ is a right coherent connected $\ZZ$-algebra, then $\grmod C$ is an abelian category, so $\tors C:=\grmod C\cap \Tors C$ is a Serre subcategory of $\grmod C$.  In this case, we define $\tails C:=\grmod C/\tors C$.

\begin{lemma} \label{lem.fito}
If $C$ is a right coherent connected $\ZZ$-algebra, then $\tors C$ is the full subcategory of $\grmod C$ consisting of finite dimensional modules.
\end{lemma}

\begin{proof}
Let $M\in \tors C$.  Since $M\in \grmod C$, there exist $x_1, \dots, x_m\in M$ such that $M=\sum _{i=1}^mx_iC$.  Since $M\in \Tors C$, $x_iC$ is right bounded for every $i=0, \dots, m$.  Since $C$ is locally finite by Lemma \ref{lem.lofi}, $x_iC$ is finite dimensional for every $i=1, \dots, m$, so $M$ is finite dimensional.  The converse is clear.
\end{proof}

\begin{remark}
If $C$ is a right coherent connected $\ZZ$-algebra, then $\tails C$ is the same as $\operatorname{cohproj} C$ defined in \cite{Na} and \cite{P} by Lemma \ref{lem.fito}, so, for the rest of the paper,
we can and will view $\tails C$ as a full subcategory of $\Tails C$ by \cite[Lemma 2.2 (2)]{Na}.
\end{remark}

We use the following $\ZZ$-algebra version of \cite[Lemma 4.3.3]{BV} (see also \cite[Proposition 1.7.11]{sheaves}) implicitly in the sequel.

\begin{lemma} \label{lem.BVl}
If $C$ is a right coherent connected $\ZZ$-algebra, then the canonical functors $\sD^b(\grmod C)\to \sD(\GrMod C)$ and $\sD^b(\tails C)\to \sD(\Tails C)$ are fully faithful.
\end{lemma}

By the above lemma, we often view  $\sD^b(\grmod C)$ and $\sD^b(\tails C)$ as full subcategories of $\sD(\GrMod C)$ and $\sD(\Tails C)$, respectively.

The next result follows from a property of a localizing subcategory.

\begin{lemma}\label{lem.tto} \textnormal {(cf. \cite [Theorem 6.8]{Na})} If $C$ is a right $\Ext$-finite connected $\ZZ$-algebra, then, for every $X\in \sD^b(\GrMod C)$, there exists a triangle
$$\R\t(X)\to X\to \R Q(X)$$
in $\sD(\GrMod C)$.  In particular, for every $M\in \GrMod C$, there exists an exact sequence
$$0\to \R^0\t (M)\to M\to \R^0Q(M)\to \R^1\t (M)\to 0,$$
and an isomorphism $\R^qQ(M)\cong \R^{q+1}\t (M)$ in $\GrMod C$ for every $q\geq 1$.
\end{lemma}

We have the following analogue of \cite[Lemma 6.9]{MN}.

\begin{lemma} \label{lemma.extendq}
Let $A, B$ be $\mathbb{Z}$-algebras.  If $B$ is a right Ext-finite connected $\mathbb{Z}$-algebra, then the left-exact functor $Q := \omega \pi:\GrMod B\to \GrMod B$ extends to a left-exact functor
$$
Q:{\Bimod }(A-B) \longrightarrow {\Bimod }(A-B)
$$
via functoriality of $Q$.
\end{lemma}

\begin{proof}
Since $\operatorname{R}^{i}\tau$ commutes with direct sums by \cite[Lemma 5.9]{MN}, then by Lemma \ref{lem.tto},  $Q$ commutes with direct sums.  Thus, if $M$ is an object in ${\Bimod }(A-B)$, then $Q(\oplus_{i}e_{i}M) \cong \oplus_{i}Q(e_{i}M)$ as graded right $B$-modules.  Furthermore, we may define a graded left $A$-module structure on this module via functorality of $Q$, and this makes $Q(M)$ an object in ${\Bimod }(A-B)$ as one can check.  It is also routine to check this defines a functor and we omit the verification.
\end{proof}

Since ${\Bimod }(A-B)$ has enough injectives, there is a right derived functor
$$
\operatorname{R}Q:{\sD}^{+}({\Bimod }(A-B)) \rightarrow {\sD}({\Bimod }(A-B)).
$$

In what follows, we shall abuse notation by writing $\tau$ and $Q$ for the associated extensions of these functors to bimodule categories.

\begin{lemma} \label{q.ww}
If $C$ is a right $\Ext$-finite connected $\ZZ$-algebra, then there exists a triangle
$$\R\t(C)\to C \to \R Q(C)$$
in $\sD(\Bimod (C - C))$.  In particular, there exist an exact sequence
$$
0\to \R^0\t(C)\to C\to \R^0Q(C)\to \R^1\t (C)\to 0,
$$
and an isomorphism
$\R^{q+1}\t(C)\cong \R^q Q(C)$
for every $q\geq1$ in $\Bimod (C - C)$.
\end{lemma}

\begin{proof}
Let $I$ denote an injective resolution of $C$ in $\Bimod (C - C)$.  For each $i \in \mathbb{Z}$, $e_{i}I$ is an injective resolution of $P_{i}$ in $\GrMod C$ by \cite[Lemma 2.3]{MN}, and thus, by Lemma \ref{lem.tto}, there is a short exact sequence of complexes
$$
0 \to \tau(e_{i}I) \to e_{i}I \to Q(e_{i}I) \to 0
$$
for each $i \in \mathbb{Z}$.  Since $\tau$ and $Q$ commute with direct sums, 
then, by definition of the extensions of $\tau$ and $Q$ to bimodules, there is a short exact sequence of complexes of objects in $\Bimod (C - C)$
$$
0 \to \tau(I) \to I \to Q(I) \to 0.
$$
The result now follows from \cite[Remark, p. 63]{hartshorne}.
\end{proof}

For a connected $\ZZ$-algebra $C$, the cohomological dimension of $\tau$ is defined by
$$\operatorname{cd }\tau:=\sup \{q\in \NN\mid \R^q\t(M)\neq 0 \textnormal { for } M\in \GrMod C\}.$$
Note that if $\operatorname{cd }\tau<\infty$, then $D\R\t(C)\in \sD^b(\Bimod (C-C))$ is bounded, which is often an essential condition in the sequel.

\begin{lemma} \label{lem.bound} Let $C$ be a right $\Ext$-finite connected $\ZZ$-algebra.  If $\operatorname{cd} \t < \infty$, then $D \R Q (C) \cong \operatorname{cone}(DC \to D \R \t (C))[-1]$ in $\sD(\Bimod (C-C))$.  In particular, $D \R Q (C)$ is bounded.
 \end{lemma}

\begin{proof} The fact that $D\R\t(C)$ is bounded follows immediately from $\operatorname{cd} \t < \infty$.
By Lemma \ref{q.ww}, there is a triangle
$$
\R\t(C)\to C\to \R Q(C)
$$
in $\sD(\operatorname{Bimod }(C-C))$.  Thus, there is a triangle
$$
D\R Q(C)\to DC\to D\R\t(C)
$$
which may be rotated to a triangle
$$
DC \to D\R \t(C) \to D \R Q(C)[1].
$$
It follows that $D \R Q (C) \cong \operatorname{cone}(DC \to D \R \t (C))[-1]$.
\end{proof}

\begin{definition} \label{def.mst}
Let $C$ be a right $\Ext$-finite connected $\ZZ$-algebra, $\cM\in \sD(\Tails C)$ and $q\in \ZZ$.
\begin{enumerate}
\item{} We define a graded right $C$-module structure on
$$
\uExt_{\cC}^q(\cC, \cM):=\oplus _{i\in \ZZ}\Hom_{\cC}(\cP_i, \cM[q])
$$

by
\begin{align*}
\uExt_{\cC}^q(\cC, \cM)_i\times C_{ij} & \cong \Hom_{\cC}(\cP_i, \cM[q])\times \Hom_C(P_{j}, P_{i}) \\
& \to \Hom_{\cC}(\cP_i, \cM[q])\times \Hom_{\cC}(\cP_j, \cP_i) \\
& \to \Hom_{\cC}(\cP_j, \cM[q]) \\
& =\uExt^q_{\cC}(\cC, \cM)_j.
\end{align*}
\item{} We define a graded left $C$-module structure on
$$
\uExt^q_{\cC}(\cM, \cC):=\oplus _{j\in \ZZ}\Hom_{\cC}(\cM[-q], \cP_j)
$$
by
\begin{align*}
C_{ij}\times \uExt^q_{\cC}(\cM, \cC)_j & \cong \Hom_C(P_j, P_i)\times \Hom_{\cC}(\cM[-q], \cP_j) \\
& \to \Hom_{\cC}(\cP_j, \cP_i)\times \Hom_{\cC}(\cM[-q], \cP_j) \\
& \to \Hom_{\cC}(\cM[-q], \cP_i) \\
& =\uExt^q_{\cC}(\cM, \cC)_i.
\end{align*}
\end{enumerate}
\end{definition}

\begin{lemma}  \label{lem.w}
If $C$ is a right $\Ext$-finite connected $\ZZ$-algebra, then $
\omega(-)\cong \uHom_{\cC}(\cC, -):\Tails C\to \GrMod C
$
as functors.  In particular, $\R^q\omega (\cM)\cong\uExt^q_{\cC}(\cC, \cM)$ in $\GrMod C$ for every $\cM\in \Tails C$ and $q\in \ZZ$.
\end{lemma}

\begin{proof} For $\cM\in \Tails C$,
\begin{align*}
\uHom_{\cC}(\cC, \cM) & := \oplus _{i\in \ZZ}\Hom_{\cC}(\cP_i, \cM)\cong \oplus _{i\in \ZZ}\Hom_C(P_i, \omega (\cM)) \\
& =: \uHom_C(C, \omega (\cM))\cong \omega (\cM)
\end{align*}
as graded vector spaces where the last isomorphism is in $\GrMod C$
by Lemma \ref{lem.Hom} (1).
Since $\pi :\GrMod C\to \Tails C$ is a functor such that $\pi \omega \cM\cong \cM$,
we have
the following commutative diagram
$$\begin{CD}
(\omega \cM)_i\times C_{ij} @>>> (\omega \cM)_j \\
\parallel & & \parallel \\
\Hom_{C}(P_i, \omega \cM)\times \Hom_C(P_j, P_i) @>>> \Hom_{C}(P_j, \omega \cM) \\
@VVV @VVV \\
\Hom_{\cC}(\cP_i, \cM)\times \Hom_{\cC}(\cP_j, \cP_i) @>>> \Hom_{\cC}(\cP_j, \cM) \\
@AAA  \parallel \\
\uHom_{\cC}(\cC, \cM)_i\times C_{ij} @>>> \uHom_{\cC}(\cC, \cM)_j \\
\end{CD}$$
so $\omega \cM\cong \uHom_{\cC}(\cC, \cM)$ in $\GrMod C$.  It follows that
$$
\omega(-)\cong \uHom_{\cC}(\cC, -):\Tails C\to \GrMod C
$$
as functors, so $\R\omega \cM\cong \RuHom_{\cC}(\cC, \cM)$ in $\sD(\GrMod C)$, hence $\R^q\omega \cM= h^q(\R\omega \cM)\cong h^q(\RuHom_{\cC}(\cC, \cM))=\uExt^q_{\cC}(\cC, \cM)$ in $\GrMod C$ for every $q\in \ZZ$.
\end{proof}

\begin{remark} \label{rem.ccbi} If $C$ is a right Ext-finite connected $\ZZ$-algebra, then $\uHom_{\cC}(\cC, \cC)$ has a structure of a bigraded $C$-bimodule as defined in Definition \ref{def.mst}.  On the other hand, since $\uHom_{\cC}(\cC, \pi (-)):\GrMod C\to \GrMod C$ is a functor commuting with direct sums, $\uHom_{\cC}(\cC, \cC)$ has a structure of a bigraded $C$-bimodule as defined in the proof of Lemma \ref{lemma.extendq}.  It is routine to check that these bigraded $C$-bimodule structures are the same.\end{remark}

The following result plays a key role in the proof of Proposition \ref{prop.asf2'}, which gives necessary and sufficient conditions for an algebra to satisfy various regularity conditions.

\begin{lemma} \label{lem.ww}
If $C$ is a right $\Ext$-finite connected $\ZZ$-algebra, then
$$
\RuHom_{\cC}(\cC, \cC) \cong \R Q(C).
$$
Therefore, there exists a triangle
$$\R\t(C)\to C\to \RuHom_{\cC}(\cC, \cC)$$
in $\sD(\Bimod (C - C))$.  In particular, there exists an exact sequence

$$
0\to \R^0\t(C)\to C\to \uHom_{\cC}(\cC, \cC)\to \R^1\t (C)\to 0
$$
and an isomorphism $\R^{q+1}\t(C)\cong \uExt^q_{\cC}(\cC, \cC)$ for every $q\geq1$ in $\Bimod (C - C)$.
\end{lemma}

\begin{proof}
The second and third part of the result will follow from the first by Lemma \ref{q.ww}.
To prove the first part of the result, we know that $Q(-)\cong \uHom_{\cC}(\cC, \pi(-)):\GrMod C\to \GrMod C$ as functors by Lemma \ref{lem.w}, so that, if we extend ${\uHom_{\cC}(\cC, \pi(-))}$ to an endofunctor on $\Bimod (C - C)$ using functoriality as we did in the proof of Lemma \ref{lemma.extendq}, we have $Q(-)\cong  {\uHom_{\cC}(\cC, \pi(-))}:\Bimod (C - C) \to \Bimod (C - C)$
as functors.
\end{proof}

\subsection{Local Duality}

We will need the following notation from \cite{MN}: viewing $k$ as a graded algebra concentrated in degree 0, we define a connected $\ZZ$-algebra $K:=\overline k$, that is,
$$
K_{ij}=k_{j-i}=\begin{cases} k & \textnormal { if } i=j \\ 0 & \textnormal { if } i\neq j.\end{cases}
$$
Note that $K^o=K$ as $\ZZ$-algebras, so $\Bimod (K^{o} - C)=\Bimod (K - C)$ denotes the category of bigraded $K$-$C$-bimodules.

In \cite{MN} and \cite{CMN}, local duality is proven for objects of $\sD^{-}(\Bimod (K-C))$.  However, for various applications in this paper, we need to apply it to objects of $\sD^{-}(\GrMod C)$.  For this reason we introduce the functor
$$
I_i(-):=Ke_i\otimes _k-:\GrMod C\to \Bimod (K- C).
$$
\begin{lemma}  \label{lem.k-c} Let $C$ be a $\ZZ$-algebra.
\begin{enumerate}
\item{} $I_i:\GrMod C\to \Bimod (K - C)$ and $e_j(-):\Bimod (K - C)\to \GrMod C$ are exact functors, which induce functors $I_i:\sD(\GrMod C)\to \sD(\Bimod (K - C))$ and $e_j(-):\sD(\Bimod (K - C))\to \sD(\GrMod C)$ such that
    $$
    e_j(-)\circ I_i\cong \begin{cases} \id & \textnormal { if } i=j \\ 0 & \textnormal { if  } i\neq j.\end{cases}
    $$

\item{}  For every $i, j\in \ZZ$, $I_i(P_j)\in \Bimod (K - C)$ is a projective bigraded $K$-$C$ bimodule.

\item{} $I_{i}$ induces a fully faithful functor $\sD(\GrMod C)\to \sD(\Bimod (K - C))$.

\end{enumerate}
\end{lemma}

\begin{proof} (1) Clearly, $I_i$ and $e_j(-)$ are exact functors.  If $M\in \GrMod C$, then
\begin{eqnarray*}
e_jI_i(M) & = & \oplus _{s\in \ZZ}(Ke_i\otimes _kM)_{js} \\
& = & \oplus _{s\in \ZZ}((Ke_i)_j\otimes _kM_s) \\
& = & \begin{cases} \oplus_{s\in \ZZ}M_s=M & \textnormal { if } i=j \\ 0 & \textnormal { if  } i\neq j.\end{cases}
\end{eqnarray*}
(2) By \cite[Lemma 2.4]{MN}, $I_i(P_j)=Ke_i\otimes _ke_jC$ is a projective bigraded $K$-$C$ bimodule.

Part (3) follows immediately from (1).
\end{proof}

Using the lemma, we have the following consequence of \cite[Theorem 2.1]{CMN}.

\begin{theorem} \label{thm.ld1}
Let $C$ be a right $\Ext$-finite connected $\ZZ$-algebra such that $\operatorname{cd }\tau < \infty$.  For $M\in \sD^{-}(\GrMod C)$,
$$D\R \t(M)\cong \RuHom_C(M, D\R \t(C))$$
in $\sD(\GrMod C^o)$.
\end{theorem}

\begin{proof}
For $M\in  \sD^{-}(\GrMod C)$,
$$\begin{array}{rll}
D\R \t(M)
& \cong D\R \t(e_0I_0(M))  & \textnormal{[Lemma \ref{lem.k-c}  (1)]} \\
& \cong D(e_0\R \t(I_0(M)))  \\
& \cong (D\R \t(I_0(M))e_0 \\
& \cong \RuHom_C(I_0(M), D\R \t(C))e_0 & \textnormal{\cite[Theorem 2.1]{CMN}} \\
& \cong \RuHom_C(e_0I_0(M), D\R \t(C)) \\
& \cong \RuHom_C(M, D\R \t(C)) & \textnormal{[Lemma \ref{lem.k-c}  (1)]}
\end{array}$$
in $\sD(\GrMod C^o)$.
\end{proof}

For a connected $\ZZ$-algebra $C$, we define the small global dimension of $C$ by
$$
\operatorname{sgldim} C:=\sup\{\pd (S_i)\mid i\in \ZZ\}.
$$
For the readers convenience, we recall \cite[Lemma 2.2]{CMN}:
\begin{lemma} \label{lemma.sgldim}
Let $C$ be a right Ext-finite connected $\mathbb{Z}$-algebra.  If $\operatorname{sgldim} C < \infty$, then $\operatorname{cd }\tau < \infty$.
\end{lemma}

\begin{lemma} \label{lemma.global}
Let $C$ be a connected $\mathbb{Z}$-algebra.  If $\operatorname {sgldim}C<\infty$, then the global dimensions of $\GrMod C$ and $\GrMod C^o$ are finite.
\end{lemma}

\begin{proof}
By \cite[Corollary 4.10]{MN}, $\operatorname{sgldim}C^o=\operatorname{sgldim}C<\infty$.
Since the functor $-\intotimesc C/C_{\geq 1} :{\GrMod }C \longrightarrow {\GrMod }C$ commutes with direct limits, the left derived functors $\uTor^{C}_{i}(-,C/C_{\geq 1})$, $i >0$ defined in \cite[Section 4.2]{MN}, also commute with direct limits by the $\mathbb{Z}$-algebra analogue of the usual argument in the graded context.  Therefore, since every graded module is a direct limit of left-bounded modules, the fact that the global dimension of $\GrMod C$ is finite follows by \cite[Proposition 4.11]{MN}, \cite[Proposition 4.7]{MN}, and the fact that $\operatorname {sgldim}C^o<\infty$.  The fact that the global dimension of $\GrMod C^o$ is finite now follows by symmetry.
\end{proof}


\begin{lemma} \label{lem.duc}
Let $C$ be a coherent connected $\ZZ$-algebra.  If $\operatorname{sgldim}C<\infty$, then we have a duality
$$\RuHom_C(-, C):\sD^b(\grmod C)\leftrightarrow \sD^b(\grmod C^o):\RuHom_{C^o}(-, C^o).$$
\end{lemma}

\begin{proof}
By \cite[Corollary 4.10]{MN}, $\operatorname{sgldim}C^o=\operatorname{sgldim}C<\infty$.  Thus, by Lemma \ref{lemma.global} and \cite[p. 68, Example 1]{hartshorne},
$$
\RuHom_{C^o}(\R\uHom_{C}(-, C), C^o): \sD({\GrMod }C) \rightarrow \sD({\GrMod }C)
$$
is way out on both sides.  Furthermore, $\uHom_{C}(-,C):{\GrMod } C \rightarrow {\GrMod }C^{o}$ equals the functor $\uHom_{C}(I_{0}(-),C)e_{0}$ when $\uHom_{C}(-,C)$ is considered as a functor from ${\Bimod }(K-C)$ to ${\Bimod }(C-K)$.   Thus,  as in the proof of \cite[Theorem 3.7]{CMN}, the functors $\RuHom_C(-, C)$ and $\RuHom_{C^o}(-, C^o)$ induce functors between the categories $\sD^b(\grmod C)$ and $\sD^b(\grmod C^o)$ by  Lemma \ref{lem.BVl}.
The result now follows from \cite[Lemma 3.6 (2)]{CMN} in the case $l=0$ and $\nu=\operatorname{id}$ (see the argument in the proof of \cite[Theorem 3.7]{CMN}). 
\end{proof}

\begin{definition} \label{def.fgf}
We say $L \in \sD^{b}({\Bimod }(K-C))$ is a perfect complex if the terms of $L$ are finite direct sums of modules of the form $Ke_{i}\otimes_{k}e_{j}C$. 
\end{definition}

\begin{lemma} \label{lem.sgl2}
Let $C$ be a right coherent connected $\ZZ$-algebra.  If $\operatorname {sgldim}C<\infty$, then every $X\in \sD^b(\grmod C)$ has a finitely generated free resolution of finite length.  Therefore, $I_{0}(X)$ is quasi-isomorphic to a perfect complex.
\end{lemma}

\begin{proof}
For $X\in \sD^b(\grmod C)$, we may assume that $X^q=0$ for all $q\gg 0$ and all $q\ll 0$.  Since $X^q\in \grmod C$ is left bounded, $\pd (X^q)<\infty$ by \cite [Proposition 4.11]{MN}, so $X^q$ has a unique minimal finitely generated free resolution of finite length for every $q\in \ZZ$ by Lemma \ref{lem.ab} and \cite[Corollary 4.5]{MN}.  The total complex of the Cartan-Eilenberg resolution of $X$ is a finitely generated free resolution of $X$ of finite length. 
\end{proof}

The following result, which will be employed to prove Theorem \ref{thm.LD3}, is a $\ZZ$-algebra version of \cite[Proposition 2.1]{J}.  The proof employs the functors $\uHom_{C}^{\bullet}(-,-)$ and $\operatorname{Tot}(- \intotimesc -)$ defined in \cite[Section 6]{MN}.

\begin{prop} \label{prop.identity}
Let $B$ and $C$ be $\mathbb{Z}$-algebras, $X \in {\sD}^{b}({\Bimod }(K-C))$, $Y \in {\sD}^{b}({\Bimod}(B-C))$ and $Z \in {\sD}^{-}({\Bimod }(K-B))$.  If $X$ is quasi-isomorphic to a perfect complex, then there is an isomorphism
$$
\operatorname{R}\uHom_C(X,Z \intotimesb^{\operatorname{L}}Y) \cong Z \intotimesb^{\operatorname{L}}\operatorname{R}\uHom_{C}(X,Y).
$$
\end{prop}

\begin{proof}
By assumption, $X$ is quasi-isomorphic to a perfect complex $L$.   Moreover, by \cite[Lemma 2.4]{MN}, $Z$ is quasi-isomorphic to a bounded above complex $F$, whose terms are direct sums of modules of the form $Ke_{i}\otimes_{k}e_{j}B$.  Thus, by \cite[Proposition 6.6 (2), Proposition 6.8]{MN},
$$
\operatorname{R}\uHom_{C}(X,Z \intotimesb^{\operatorname{L}}Y) \cong \uHom_{C}^{\bullet}(L,\operatorname{Tot}(F \intotimesb Y)),
$$
and
$$
Z \intotimesb^{\operatorname{L}}\operatorname{R}\uHom_{C}(X,Y) \cong \operatorname{Tot}(F \intotimesb \uHom_{C}^{\bullet}(L,Y)).
$$
Therefore, it suffices to prove that there is a natural isomorphism of complexes
$$
\uHom_{C}^{\bullet}(L,\operatorname{Tot}(F \intotimesb Y)) \cong \operatorname{Tot}(F \intotimesb \uHom_{C}^{\bullet}(L,Y)).
$$
Since, for $P$ a free module in ${\GrMod }B$, $M$ in ${\GrMod }C$ a finitely generated free module, and $N$ in ${\Bimod }(B-C)$, there is a natural isomorphism
$$
 P\intotimesb \uHom_{C}(M,N) \longrightarrow \uHom_{C}(M, P \intotimesb N),
$$
in ${\GrMod }C$, the remainder of the proof is the same as that of \cite[Proposition 2.1]{J}.
\end{proof}


We also will employ the following variant of \cite[Lemma 6.10]{MN}. 

\begin{lemma} \label{lemma.ldq}
Let $C$ be a right Ext-finite connected $\mathbb{Z}$-algebra.  If $M$ is a complex of free right $C$-modules, and $N$ is a complex of bigraded $C$-bimodules, then there exists a canonical isomorphism
$$
\operatorname{Tot}(M \intotimesc Q(N)) \cong Q(\operatorname{Tot}(M
 \intotimesc N)).
$$
\end{lemma}

\begin{proof}
By Lemma \ref{lem.tto} and \cite[Lemma 5.9]{MN},
$\R^{i}Q$ commutes with direct limits for $i \geq 0$.  Thus, the argument in \cite[Lemma 6.10]{MN} can be applied to prove the result.
\end{proof}

The following is a $\mathbb{Z}$-algebra version of local duality found in \cite{NV}.

\begin{theorem} \label{thm.ld0}
Let $C$ be a right $\Ext$-finite connected $\ZZ$-algebra such that $\operatorname{cd }\tau<\infty$.
For $M \in \sD^{-}(\Bimod (K-C))$, 
$$D\R Q(M)\cong \RuHom_C(M, D\R Q(C))$$
in $\sD(\Bimod(C-K))$.
\end{theorem}

\begin{proof}
Since $\operatorname{cd }\tau < \infty$, the cohomological dimension of $Q$ is finite by Lemma \ref{lem.tto}.  Furthermore, by Lemma \ref{lem.tto} and \cite[Lemma 5.9]{MN}, $\R^{i}Q$ commutes with direct limits for $i \geq 0$.  Thus, by Lemma \ref{lemma.ldq}, the argument of \cite[Theorem 2.1]{CMN} can be applied to prove the result.
\end{proof}

\section{AS-regular $\ZZ$-algebras} \label{section.regularity}

AS-regular algebras were originally introduced in \cite{AS}, and play an essential role in noncommutative algebraic geometry.  The related notion of ASF-regular algebras was originally introduced in \cite{MM}.  After recalling these notions, we modify these definitions for the purpose of this paper.

\begin{definition} \label{def.ASg}
Let $A$ be a locally finite connected graded algebra.
\begin{enumerate}
\item{} $A$ is called {\it AS-regular of dimension $d$ and of Gorenstein parameter $\ell$} if
\begin{itemize}
\item{} $\gldim A=d$, and
\item{} $\Ext^q_A(k, A(j))\cong \begin{cases}
k & \textnormal { if $q=d$ and $j=-\ell$,} \\
0 & \textnormal { otherwise.}
\end{cases}$
\end{itemize}
\item{} $A$ is called {\it ASF-regular of dimension $d$ and of Gorenstein parameter $\ell$} if
\begin{itemize}
\item{} $\gldim A=d$, and
\item{} $D\R^q\t(A)\cong \begin{cases}
A(-\ell) & \textnormal { if $q=d$,} \\
0 & \textnormal { otherwise}
\end{cases}$
as graded left and right $A$-modules.
\end{itemize}
\end{enumerate}
\end{definition}

Recall that a connected graded algebra $A$ is right $\Ext$-finite if and only if $A$ is left $\Ext$-finite
by \cite{Ve}, so we may just say that $A$ is $\Ext$-finite, and, in this case, $A$ is locally finite.  If $A$ is an $\Ext$-finite connected graded algebra, then $A$ is AS-regular of dimension $d$ and of Gorenstein parameter $\ell$ if and only if  $A$ is ASF-regular of dimension $d$ and of Gorenstein parameter $\ell$
by \cite[Theorem 3.12]{MM} and \cite[Theorem 2.19]{MU1}.

In this section, we extend these definitions to $\ZZ$-algebras, and study some properties and relationships of these $\ZZ$-algebras.

\subsection{AS-regular $\ZZ$-algebras and ASF-regular $\ZZ$-algebras}
The following definition was given in \cite{MN}.

\begin{definition} \textnormal{\cite[Definition 7.1]{MN}} A locally finite connected $\ZZ$-algebra $C$ is called {\it AS-regular of dimension $d$ and of Gorenstein parameter $\ell$} if
\begin{enumerate}
\item[(ASR1)] $\pd_C S_i=d$ for every $i\in \ZZ$, and
\item[(ASR2)] $\Ext^q_C(S_i, P_j)=\begin{cases} k & \textnormal { if $q=d$ and $j=i+\ell$} \\
0 & \textnormal { otherwise,} \end{cases}$
that is,
$$
\RHom_C(S_i, P_j)\cong \begin{cases} k[-d] & \textnormal { if } j=i+\ell \\
0 & \textnormal {otherwise} \end{cases}
$$
in $\sD(\Mod k)$.
\end{enumerate}

\end{definition}

\begin{remark} \label{rem.asrzr}
Since $\uExt^q_C(S_i, C)_j=\Ext^q_C(S_i, e_jC)=\Ext^q_C(S_i, P_j)$, the condition (ASR2) is equivalent to
$$\uExt^q_C(S_i, C)\cong \begin{cases} S_{i+\ell} & \textnormal { if $q=d$} \\
0 & \textnormal { otherwise} \end{cases}$$
as graded left $C$-modules for every $i\in \ZZ$.
\end{remark}
The following proposition justifies our definition of an AS-regular $\ZZ$-algebra above. A similar result was stated in \cite[Section 4.1]{Vq} although the definition of an AS-regular $\ZZ$-algebra in \cite[Section 4.1]{Vq}  differs slightly from ours. (In \cite{Vq}, the Gorenstein parameter is not well-defined.)

\begin{proposition}
\label{prop.Vq}
Let $A$ be a locally finite connected graded algebra.  Then $A$ is an AS-regular algebra
of dimension $d$ and of Gorenstein parameter $\ell$
if and only if $\overline A$ is an AS-regular $\ZZ$-algebra
of dimension $d$ and of Gorenstein parameter $\ell$.
\end{proposition}

It follows from the above proposition and Lemma \ref{lem.1p} that $C$ is a 1-periodic AS-regular $\ZZ$-algebra of dimension $d$ and of Gorenstein parameter $\ell$ if and only if
there exists an AS-regular algebra $A$ of dimension $d$ and of Gorenstein parameter $\ell$ such that $C\cong \overline A$.

The next result tells us that AS-regularlity is left-right symmetric.

\begin{theorem} \label{thm.lras}
Let $C$ be a 
connected $\ZZ$-algebra.  Then $C$ is a right $\Ext$-finite AS-regular algebra of dimension $d$ and of Gorenstein parameter $\ell$ if and only if $C^{o}$ is a right $\Ext$-finite AS-regular algebra of dimension $d$ and of Gorenstein parameter $\ell$.
\end{theorem}

\begin{proof}
Let $C$ be a right $\Ext$-finite AS-regular $\ZZ$-algebra of dimension $d$ and of Gorenstein parameter $\ell$.
Since $C$ is right $\Ext$-finite, $S_i$ has a minimal finitely generated free resolution
$$0\to F^d\to \cdots \to F^0\to S_i\to 0$$
in $\GrMod C$.
Since $F^q\in \add \{P_{j}\}_{j \in \mathbb{Z}}$, we have $\uHom_C(F^q, C)\in \add\{Q_j\}_{j\in \ZZ}=\add\{P_j^o\}_{j\in \ZZ}$ for every $q$ by Lemma \ref{lem.Hom} (2).  By Remark \ref{rem.asrzr},
$$0\to \uHom_C(F^0, C)
\to \cdots \to \uHom_C(F^d, C)\to \uExt^d_C(S_i, C)\cong S_{i+\ell}\to 0$$ is a minimal finitely generated free resolution of $S_{-\ell-i}^{o}\cong S_{i+\ell}$ in $\GrMod C^o$, so 
$C^o$ is right $\Ext$-finite and 
$\pd S_{-\ell-i}^o=d$ for every $i\in \ZZ$.

Since $\R \uHom_{C^o}(\R \uHom_C(F, C), C^o)\cong F$ by \cite[Lemma 3.6 (2)]{CMN} (see Lemma \ref{lem.duc}),
$$0\to \uHom_{C^o}(\uHom_C(F^d, C), C^o)\to \cdots \to \uHom_{C^o}(\uHom_C(F^0, C), C^o)\to S_{i}\to 0$$
is isomorphic to the original minimal finitely generated free resolution of $S_i$ in $\GrMod C$ by the uniqueness of the minimal free resolution, so

$$\uExt^q_{C^{o}}(S^{o}_{-\ell-i}, C^o)=\begin{cases} S_i\cong S^o_{-i} & \textnormal { if $q=d$} \\
0 & \textnormal { otherwise,} \end{cases}$$
hence $C^{o}$ is an AS-regular $\ZZ$-algebra of dimension $d$ and of Gorenstein parameter $\ell$ by Remark \ref{rem.asrzr}.
\end{proof}

\begin{remark} \label{rem.lras} In the setting of the above theorem, since
$$
\Hom_C(P_i, S_j)=\begin{cases} k & \textnormal { if } i=j \\
0 & \textnormal { if } i\neq j\end{cases}
$$
by Lemma \ref{lem.Hom} (1), we have $F^0\cong P_i$.  Similarly, since
$$
\Hom_{C^o}(P_i^o, S_j^o)=\begin{cases} k & \textnormal { if } i=j \\
0 & \textnormal { if } i\neq j,\end{cases}
$$
we have $\uHom_C(F^d, C)\cong P_{-\ell-i}^o\cong Q_{i+\ell}$, so we have $F^d\cong P_{i+\ell}$.  It follows that the minimal finitely generated free resolution of $S_i$ in $\GrMod C$ is of the form
$$0\to P_{i+\ell}\to \cdots \to P_i\to S_i\to 0.$$
\end{remark}

The following definition was also given in \cite{MN}.

\begin{definition} \textnormal{\cite[Definition 7.5]{MN}} \label{def.ASF}
A connected $\ZZ$-algebra $C$ is called {\it ASF-regular of dimension $d$ and of Gorenstein parameter $\ell$} if
\begin{enumerate}
\item[(ASF1)] $\operatorname {sgldim} C
=d<\infty$, and
\item[(ASF2)] $\R^q\t(P_j)\cong \begin{cases}
D(Q_{j-\ell}) & \textnormal {  if $q=d$} \\
0 & \textnormal { otherwise}
\end{cases}$ as graded right $C$-modules for every $j\in \ZZ$, that is,
$$
\R\t(P_j)\cong (D(Q_{j-\ell}))[-d]
$$
in $\sD(\GrMod C)$.
\end{enumerate}
\end{definition}

As to relationships between AS-regular $\ZZ$-algebras and ASF-regular $\ZZ$-algebras, one implication was proved in \cite{MN}.

\begin{theorem}  \textnormal{\cite[Corollary 7.7]{MN}}  \label{thm.asasf}
If $C$ is a right
$\Ext$-finite AS-regular $\ZZ$-algebra of dimension $d$ and of Gorenstein parameter $\ell$, then $C$ is an ASF-regular $\ZZ$-algebra of dimension $d$ and of Gorenstein parameter $\ell$.
\end{theorem}

\subsection{ASF$^+$-regular $\ZZ$-algebras}
We now modify the original definition of an ASF-regular algebra given in Definition \ref{def.ASg}, replacing the condition (ASF2) with conditions that takes into account both the left and right $C$-module structure on $\R\t(C)$.

The following definition is closer to the original definition of an ASF-regular algebra.

\begin{definition} \label{def.asf+} A locally finite connected $\ZZ$-algebra $C$ is called {\it ASF$^+$-regular of dimension $d$ and of Gorenstein parameter $\ell$} if
\begin{enumerate}
\item[(ASF1)] $\operatorname {sgldim} C=d<\infty$,

\item[(ASF2)]
$(D\R\t(C))e_j\cong Ce_{j-\ell}[d]=Q_{j-\ell}[d]$ in $\sD(\GrMod C^o)$ for every $j\in \ZZ$, and
\item[(ASF2$^+$)]
$e_i(D\R\t(C))\cong e_{i+\ell}C[d]=P_{i+\ell}[d]$
in $\sD(\GrMod C)$ for every $i\in \ZZ$.
\end{enumerate}
\end{definition}

\begin{remark} \label{rem.lrl}
Let $C$ be a connected $\ZZ$-algebra.
\begin{enumerate}
\item{} Since
$$(D\R\t(C))e_j\cong D(e_j\R\t(C))\cong D(\R\t(e_jC))=D\R\t(P_j)$$
as graded left $C$-modules, the above condition (ASF2) is equivalent to the condition (ASF2) in Defintion \ref{def.ASF}.  However, since $Q_j$ has no graded right $C$-module structure, $\t(Q_j)$ is not well-defined, so we are not able to replace $e_i(D\R\t(C))$ by $D\R\t(Q_i)$ in the above condition (ASF2$^+$).

\item{}
Since $C(0, -\ell)e_j=\oplus _{i\in \ZZ}C_{i, j-\ell}=Q_{j-\ell}$ as graded left $C$-modules and  $e_iC(\ell,0)=\oplus _{j\in \ZZ}C_{i+\ell, j}=P_{i+\ell}$ as graded right $C$-modules,
the condition (ASF2) is equivalent to
$$
D\R^q\t(C) \cong \begin{cases} C(0,-\ell) & \textnormal{ if }  q=d \\
0 & \textnormal { if } q\neq d \end{cases}
$$
as graded left $C$-modules, and the condition (ASF2$^+$)  is equivalent to
$$
D\R^q\t(C) \cong \begin{cases} C(\ell,0) & \textnormal{ if }  q=d \\
0 & \textnormal { if } q\neq d \end{cases}
$$
as graded right $C$-modules.

\item{}

Recall that $C(0,-\ell)$ has a bigraded $C$-$C(-\ell)$ bimodule structure, and $C(\ell,0)$ has a bigraded $C(\ell)$-$C$ bimodule structure.  Although $D\R^d\t(C)$ has a bigraded $C$-bimodule structure, $C(0,-\ell), C(\ell,0)$ do not have natural bigraded $C$-bimodule structures unless $C$ is $\ell$-periodic.  Even if $C$ is $\ell$-periodic, $C(0,-\ell)$ and $C(\ell,0)$ are not isomorphic to $C$ as bigraded $C$-bimodules in general.
\end{enumerate}
\end{remark}


\begin{proposition} \label{prop.asf2'}
Let $C$ be a right $\Ext$-finite
connected $\ZZ$-algebra of $\operatorname{sgldim} C=d$.
\begin{enumerate}
\item{} For $d=1$,
\begin{itemize}
\item{} $C$ satisfies (ASF2) if and only if there is an exact sequence
$$0\to C\to \uHom_{\cC}(\cC, \cC)\to D(C(0,-\ell))\to 0$$
as graded right $C$-modules and $\uExt^q_{\cC}(\cC, \cC)=0$ for every $q\geq 1$.
\item{} $C$ satisfies (ASF2$^+$) if and only if there is an exact sequence
$$0\to C\to \uHom_{\cC}(\cC, \cC)\to D(C(\ell,0))\to 0$$
as graded left $C$-modules and $\uExt^q_{\cC}(\cC, \cC)=0$ for every $q\geq 1$.
\end{itemize}
\item{} For $d\geq 2$,
\begin{itemize}
\item{} $C$ satisfies (ASF2) if and only if
$$\uExt^q_{\cC}(\cC, \cC)\cong \begin{cases} C & \textnormal { if $q=0$} \\
D(C(0,-\ell)) & \textnormal { if $q=d-1$} \\
0 & \textnormal { otherwise}
\end{cases}$$ as graded right $C$-modules, and
\item{} $C$ satisfies (ASF2$^+$) if and only if
$$
\uExt^q_{\cC}(\cC, \cC)\cong \begin{cases} C & \textnormal { if $q=0$} \\
D(C(\ell,0)) & \textnormal { if $q=d-1$} \\
0 & \textnormal { otherwise}
\end{cases}$$ as graded left $C$-modules.
\end{itemize}
\end{enumerate}
In the above, we tacitely require that the isomorphism when $q=0$ is induced by the canonical map
$C\cong \uHom_C(C, C)\to \uHom_{\cC}(\cC, \cC)$.
\end{proposition}

\begin{proof}
By Lemma \ref{lem.ww} and Lemma \ref{rem.lrl} (2), we have a triangle
$$(D(C(0,-\ell)))[-d]\cong \R\t(C)\to C\to \R\uHom_{\cC}(\cC, \cC)$$
in $\sD(\GrMod C)$ and a triangle
$$(D(C(\ell,0)))[-d]\cong \R\t(C)\to C\to \R\uHom_{\cC}(\cC, \cC)$$
in $\sD(\GrMod C^o)$, so the result follows.
\end{proof}

\begin{lemma} \label{lem.asfz2}
Let $A$ be an $\Ext$-finite connected graded algebra, and $\Phi:\GrMod A\to \GrMod \overline A, \Phi^o:\GrMod A^o\to \GrMod \overline A^o$ equivalence functors defined in Lemma \ref{lem.ovwi} and Lemma \ref{lem.asfz1}.
\begin{enumerate}
\item{} $\overline A$ is right $\Ext$-finite.
\item{} $\Phi\t\cong \t \Phi:\GrMod A\to \GrMod \overline A$ as functors.

\item{} For $M\in \GrMod A^e$, $e_i{\R\t(\overline M)}\cong \Phi(\R\t(M(-i)))$ for every $i\in \ZZ$ in $\sD(\GrMod \overline A)$, and ${\R\t(\overline M)}e_j\cong \Phi^o(\R\t(M(j)))$ for every $j\in \ZZ$ in $\sD(\GrMod \overline A^o)$.
\end{enumerate}
\end{lemma}

\begin{proof}  (1) Since $\Phi(A(-i))\cong P_i$, if $F$ is a minimal finitely generated free resoution of $k(j)$ in $\GrMod A$, then $\Phi(F)$ is a minimal finitely generated free resoution of $\Phi(k(-i))\cong S_i$ in $\GrMod \overline A$ by Lemma \ref{lem.ovwi}.

(2)
Since $M$ is right bounded if and only if $\Phi(M)$ is right bounded, $\Phi$ restricts to an equivalence functor $\Phi:\Tors A\to \Tors \overline A$.
Since $\t:\GrMod A\to \Tors A,$ and $\t:\GrMod \overline A\to \Tors \overline A$ are right adjoint to the inclusion functors $\Tors A\to \GrMod A$ and $\Tors \overline A\to \GrMod \overline A$, respectively, we have $\Phi\t\cong \t \Phi:\GrMod A\to \GrMod \overline A$ as functors.

(3) For $M\in \GrMod A^e$,
$$e_i\t(\overline M)\cong \t(e_i\overline M)\cong \t \Phi(M(-i))\cong \Phi\t(M(-i))$$
in $\GrMod \overline A$ by Lemma \ref{lem.asfz1} (2), and  (2) above, so
$$e_i(-)\circ \t\circ \overline {(-)}\cong \Phi\circ \t\circ (-i):\GrMod A^e\to \GrMod \overline A$$
as functors for every $i\in \ZZ$.
Since functors $\overline {(-)}: \GrMod A^e\to \Bimod (\overline A - \overline A)$, $e_i(-):\GrMod \overline A^e\to \GrMod \overline A$, $(-i):\GrMod A^e\to \GrMod A^e$, and $\Phi:\GrMod A^e\to \GrMod \overline A$ are all exact, $e_i{\R\t(\overline M)}\cong \Phi(\R\t(M(-i)))$ for every $i\in \ZZ$ in $\sD(\GrMod \overline A)$.

On the other hand, since $\overline {(-)}:\GrMod A\to \GrMod \overline A$ is an exact functor, and
$$
(\overline {A_{\geq n}})_{ij}=
(A_{\geq n})_{j-i}=
\left\{\begin{matrix} A_{j-i} & \textnormal { if } j-i\geq n
\\
0 & \textnormal { if } j-i<n
\end{matrix}\right\}=\left\{\begin{matrix} \overline A_{ij} & \textnormal { if } j-i\geq n
\\
0 & \textnormal { if } j-i<n
\end{matrix}\right\}=(\overline A_{\geq n})_{ij}$$
for every $n\in \ZZ$,
$$\overline {A/A_{\geq n}}\cong \overline A/\overline {A_{\geq n}}\cong \overline A/\overline A_{\geq n},$$
so
\begin{eqnarray*}
\t(\overline M)e_j & \cong & \lim_{n\to \infty}  \uHom_{\overline A}(\overline A/\overline {A}_{\geq n}, \overline M)e_j \\
& = & \lim_{n\to \infty} \uHom_{\overline A}(e_j(\overline A/\overline {A}_{\geq n}), \overline M)  \\
& := & \lim_{n\to \infty}  (\oplus _{i\in \ZZ}\Hom_{\overline A}(e_j\overline {A/A_{\geq n}}, e_i\overline M)) \\
& \cong & \lim_{n\to \infty}  (\oplus _{i\in \ZZ}\Hom_{\overline A}(\Phi((A/A_{\geq n})(-j)), \Phi(M(-i))) \\
& \cong & \lim_{n\to \infty} \Phi^o(\oplus _{i\in \ZZ}\Hom_A((A/A_{\geq n})(-j), M(i))) \\
& = & \Phi^o(\lim _{n\to \infty}\uHom_A(A/A_{\geq n}, M(j))) \\
& \cong & \Phi^o(\t(M(j)))
\end{eqnarray*}
in $\GrMod \overline A^o$ by \cite[Lemma 5.8]{MN}, Lemma \ref{lem.asfz1} (2), and Lemma \ref{lem.ovwi},
hence
$$(-)e_j\circ \t\circ \overline {(-)}\cong \Phi^o\circ \t\circ (j):\GrMod A^e\to \GrMod \overline A^o$$
as functors for every $j\in \ZZ$.  

Since the functors $\overline {(-)}: \GrMod A^e\to \Bimod (\overline A - \overline A)$,
$(-)e_j:\GrMod \overline A^e\to \GrMod \overline A^o$,
$(j):\GrMod A^e\to \GrMod A^e$,
and $\Phi^o:\GrMod A^e\to \GrMod \overline A^o$
are all exact, ${\R\t(\overline M)}e_j\cong \Phi^o(\R\t(M(j)))$ for every $i\in \ZZ$ in $\sD(\GrMod \overline A^o)$.
\end{proof}

The following theorem justifies the definition of an ASF$^+$-regular $\ZZ$-algebra.

\begin{theorem} \label{thm.asfo}
Let $A$ be an $\Ext$-finite connected graded algebra.  Then $A$ is an ASF-regular algebra
of dimension $d\geq 1$ and of Gorenstein parameter $\ell$
if and only if $\overline A$ is an ASF$^+$-regular $\ZZ$-algebra
of dimension $d\geq 1$ and of Gorenstein parameter $\ell$.
\end{theorem}

\begin{proof}  If $\Phi:\GrMod \overline A\cong \GrMod A$ is an equivalence functor defined in Lemma \ref{lem.ovwi}, then $\Phi(k(-i))\cong S_i$, so
$$\pd_{\overline A}S_i=\pd _Ak(-i)=\gldim A$$
for every $i\in \ZZ$.

Since $\Phi^o(A(j-\ell))\cong \overline Ae_{j-\ell}$ in $\GrMod \overline A^o$ by Lemma \ref{lem.asfz1} (1), and
\begin{eqnarray*}
(D\R\t(\overline A))e_j & \cong & D(e_j\R\t(\overline A))\\
& \cong & D\Phi(\R\t(A(-j))) \\
& \cong & \Phi^oD\R\t(A(-j))
\end{eqnarray*}
in $\sD(\GrMod \overline A^o)$ by Lemma \ref{lem.asfz2} (3) and Lemma \ref{lem.asfz1} (3),
$$
(D\R\t(\overline A))e_j\cong \overline Ae_{j-\ell}[d]
$$
for every $j\in \ZZ$ if and only if
$$
D\R\t(A(-j))\cong A(j-\ell)[d]
$$
in $\sD(\GrMod A^o)$ for every $j\in \ZZ$ by Lemma \ref{lem.asfz1}  (1) if and only if  $$
D\R\t(A)\cong A(-\ell)[d]
$$
in $\sD(\GrMod A^o)$.

Similarly, since $\Phi(A(-i-\ell))\cong e_{i+\ell}\overline A$ in $\GrMod \overline A$ by Lemma \ref{lem.ovwi}, and
$$e_i(D\R\t(\overline A))\cong D(\R\t(\overline A)e_i)\cong D\Phi^o(\R\t(A(i)))\cong \Phi D\R\t(A(i))$$
in $\sD(\GrMod \overline A)$ by Lemma \ref{lem.asfz2} (3) and Lemma \ref{lem.asfz1} (3), $$
e_i(D\R\t(\overline A))\cong e_{i+\ell}\overline A[d]
$$
for every $i\in \ZZ$ if and only if $D\R\t(A(i))\cong A(-i-\ell)[d]$ in $\sD(\GrMod A)$ for every $i\in \ZZ$ by Lemma \ref{lem.ovwi} if and only if $D\R\t(A)\cong A(-\ell)[d]$ in $\sD(\GrMod A)$.
\end{proof}

We have the following implication.

\begin{theorem} \label{thm.asfas2}
If $C$ is a
right $\Ext$-finite ASF$^+$-regular $\ZZ$-algebra of dimension $d$ and of Gorenstein parameter $\ell$,
then $C$ is an AS-regular $\ZZ$-algebra of dimension $d$ and of Gorenstein parameter $\ell$.
\end{theorem}

\begin{proof}
First note that in this situation, $\operatorname{cd }\tau < \infty$ by Lemma \ref{lemma.sgldim}.  Since
\begin{align*}
\R\Hom_{C}(S_{i}, P_{j})
& \cong \R\Hom_{C}(S_{i}, e_{j-\ell}(D\R\t(C))[-d]) \\
& \cong \R\uHom_{C}(S_{i}, D\R\t(C))_{j-\ell}[-d] \\
& \cong D\R\t(S_{i})_{j-\ell}[-d] \\
& \cong (DS_{i})_{j-\ell}[-d] \\
& \cong \begin{cases} k[-d] & \textnormal { if } j=i+\ell \\
0 & \textnormal {otherwise} \end{cases}
\end{align*}
in $\sD(\Mod k)$ by (ASF2$^+$) and Theorem \ref{thm.ld1}, $C$ satisfies (ASR2).  Since $\pd S_i\leq \operatorname {sgldim}C=d$ by (ASF1), while $\pd S_i\geq d$ by (ASR2) for every $i\in \ZZ$, $C$ satisfies (ASR1).
\end{proof}

\subsection{ASF$^{++}$-regular $\ZZ$-algebras}

We introduce another notion of regularity, which plays an essential role in this paper.

\begin{definition} \label{def.asf++} A locally finite connected $\ZZ$-algebra $C$ is called ASF$^{++}$-regular of dimension $d$ and of Gorenstein parameter $\ell$ if
\begin{enumerate}
\item[(ASF1)] $\operatorname {sgldim} C=d<\infty$,  and
\item[(ASF2$^{++}$)] $D\R\t(C)\cong C(0,-\ell)_{\nu}[d]$ in $\sD(\Bimod (C - C))$ for some isomorphism $\nu:C\to C(-\ell)$ of $\ZZ$-algebras, called the Nakayama isomorphism of $C$.
\end{enumerate}
\end{definition}

\begin{remark}  Let $C$ be a $\ZZ$-algebra and $\l=\{\l_i\in C_{ii}\}_{i\in \ZZ}$ is a sequence of unit elements, then $I_{\l}:C\to C$ defined by $I_{\l}(a)=\l_i^{-1}a\l_j$ for $a\in C_{ij}$ is an automorphism of a $\ZZ$-algebra $C$, called an inner automorphism of $C$.  If
$$
\nu, \nu':C\to C(-\ell)
$$
are isomorphisms of $\ZZ$-algebras, and
$$
\varphi:C(0,-\ell)_{\nu}\to C(0,-\ell)_{\nu'}
$$
is an isomorphism of bigraded $C$-bimodules, then
\begin{align*}
\varphi(e_{i-\ell})\nu'(a) & = \varphi(e_{i-\ell})*a = \varphi (e_{i-\ell}*a) = \varphi (e_{i-\ell}\nu(a)) \\
& = \varphi (\nu(a)) = \varphi (\nu(a)e_{j-\ell}) = \nu(a)\varphi (e_{j-\ell})
\end{align*}
for $a\in (C(0,-\ell)_{\nu})_{i, j+\ell}=C_{ij}$ so that $\nu(a), \nu'(a)\in C(-\ell)_{ij}=C_{i-\ell, j-\ell}$, so $\nu'=I_{\l}\circ \nu$ where $\l=\{\varphi(e_{i-\ell})\in C(-\ell)_{ii}\}_{i\in \ZZ}$.
It follows that the Nakayama isomorphism is unique up to inner automorphism of $C(-\ell)$
if it exists.
\end{remark}

\begin{remark} If $A$ is an AS-regular algebra, then $D\R\t(A)\cong A(-\ell)_{\nu}[d]$ in $\sD(\GrMod A^e)$ for some graded algebra automorphism $\nu\in \Aut A$ called the Nakayama automorphism of $A$.  We say that $A$ is Calabi-Yau if $\nu=\id$ (up to inner automorphism).  For a non-Calabi-Yau AS-regular algebra $A$, it is often the case that there exists
a Calabi-Yau AS-regular algebra $A'$ such that $\GrMod A\cong \GrMod A'$ (see \cite[Theorem 1.1]{IM}).
In this case, $\overline A\cong \overline {A'}$ as $\ZZ$-algebras by \cite[Theorem 1.1]{S}, so the notion of Calabi-Yau does not make sense for an AS-regular $\ZZ$-algebra.
 (The Nakayama isomorphism for a $\ZZ$-algebra is never the identity.)
\end{remark}

\begin{lemma} \label{lem.twi}
Let $C$ be a locally finite connected $\ell$-periodic $\ZZ$-algebra,
and $M$ a bigraded $C$-bimodule.  If $M\cong C(0,-\ell)$
as graded left $C$-modules and $M\cong C(\ell,0)$
as graded right $C$-modules, then there exists an isomorphism of $\ZZ$-algebras $\nu:C\to C(-\ell)$ such that $M\cong C(0,-\ell)_{\nu}$ as bigraded $C$-bimodules.
\end{lemma}

\begin{proof}  Let $\psi:M\to C(0,-\ell)$ be an isomorphism of graded left $C$-modules.  Define a bigraded $k$-linear map $\nu:C\to C(-\ell)$ by $\nu(a)=\psi(\psi^{-1}(e_{i-\ell})a)$ for $a\in C_{ij}$.
For $a\in C_{ij}, b\in C_{jk}$,
\begin{align*}
\nu(a)\nu(b) & =\psi(\psi^{-1}(e_{i-\ell})a)\psi(\psi^{-1}(e_{j-\ell})b) \\
& =\psi(\psi(\psi^{-1}(e_{i-\ell})a)\psi^{-1}(e_{j-\ell})b) \\
& =\psi(\psi^{-1}(\psi(\psi^{-1}(e_{i-\ell})a))b) \\
& =\psi(\psi^{-1}(e_{i-\ell})ab)=\nu(ab),
\end{align*}
so $\nu:C\to C(-\ell)$ is a homomorphism of $\ZZ$-algebras.

Since $C$ is $\ell$-periodic, there exists an isomorphism $\phi:C\to C(-\ell)$ of $\ZZ$-algebras.  By Lemma \ref{lem.rel} (3), there exists an isomorphism
$$
\varphi :M\to C(\ell,0)\to C(0, -\ell)_{\phi}
$$
of graded right $C$-modules.  If $\nu(a)=0$, then $\psi^{-1}(e_{i-\ell})a=0$, so
$$\varphi(\psi^{-1}(e_{i-\ell}))\phi(a)=\varphi(\psi^{-1}(e_{i-\ell}))*a=
\varphi(\psi^{-1}(e_{i-\ell})a)=0.$$
Since $0\neq \varphi(\psi^{-1}(e_{i-\ell}))\in C_{i-\ell}\cong k$, $\phi(a)=0$, so $a=0$, hence $\nu$ is injective.  Since $C$ is locally finite, $\nu$ is surjective, so $\nu :C\to C(-\ell)$ is an isomorphism of $\ZZ$-algebras.

We now consider the bigraded $k$-linear map $\psi:M\to C(0,-\ell)_{\nu}$. For $a\in C_{ij}, m\in M_{jt}, b\in C_{ts}$,
\begin{align*}
a\psi(m)*b
& =a\psi(m)\nu(b) =\psi(am)\psi(\psi^{-1}(e_{t-\ell})b) \\
& = \psi (\psi(am)\psi^{-1}(e_{t-\ell})b) = \psi(\psi^{-1}(\psi(am))b) \\
& = \psi(amb),
\end{align*}
so $\psi:M\to C(0,-\ell)_{\nu}$ is an isomorphism of bigraded $C$-bimodules.
\end{proof}

\begin{theorem} \label{thm.++}
A locally finite connected $\ZZ$-algebra is ASF$^{++}$-regular of dimension $d$ and of Gorenstein parameter $\ell$ if and only if it is $\ell$-periodic ASF$^+$-regular of dimension $d$ and of Gorenstein parameter $\ell$.
\end{theorem}

\begin{proof} If $C$ is an ASF$^{++}$-regular $\ZZ$-algebra of dimension $d$ and of Gorenstein parameter $\ell$ with the Nakayama isomorphism $\nu:C\to C(-\ell)$, then $C$ is clearly $\ell$-periodic.  Since
$$
D\R\t(C)\cong C(0,-\ell)_{\nu}[d]
$$
in $\sD(\Bimod (C - C))$,
$$
D\R\t(C)e_j\cong C(0,-\ell)e_j[d]\cong Q_{j-\ell}[d]
$$
in $\sD(\GrMod C^o)$.  Similarly, since
$$
D\R\t(C)\cong C(0,-\ell)_{\nu}[d]\cong {_{\nu^{-1}}C(\ell,0)}[d]
$$ in $\sD(\Bimod (C - C))$ by Lemma \ref{lem.rel} (3),
$$
e_iD\R\t(C)\cong e_iC(\ell,0)[d]\cong P_{i+\ell}[d]
$$
in $\sD(\GrMod C)$, so $C$ is ASF$^+$-regular.

Conversely, let $C$ be an $\ell$-periodic ASF$^{+}$-regular $\ZZ$-algebra of dimension $d$ and of Gorenstein parameter $\ell$.  By Remark \ref{rem.lrl} (2), $D\R^d\t(C)\cong C(0,-\ell)$
as graded left $C$-modules and $D\R^d\t(C)\cong C(\ell,0)$
as graded right $C$-modules.   By Lemma \ref{lem.twi},
$$
D\R^q\t(C)\cong \begin{cases} C(0,-\ell)_{\nu} & \textnormal{ if } q=d \\
0 & \textnormal{ if } q\neq d \end{cases}
$$
for some isomorphism $\nu:C\to C(-\ell)$, so $C$ is ASF$^{++}$-regular.
\end{proof}

In summary, for a right $\Ext$-finite
connected $\ZZ$-algebra $C$, we have the following implications by Theorem \ref{thm.asasf}, Theorem \ref{thm.asfas2}, and Theorem \ref{thm.++}:
\begin{center}
\begin{tabular}{cccccccccccc}
ASF$^{++}$ & $\Rightarrow$ & ASF$^+$ & $\Rightarrow$ & AS & $\Rightarrow$ & ASF \\
\end{tabular}
\end{center}
Moreover, if $C$ has a ``balanced dualizing complex", then AS $\Leftrightarrow $ ASF by \cite[Theorem 7.10]{MN},
and if $C$ is $\ell$-periodic, then ASF$^{++}$ $\Leftrightarrow$ ASF$^+$ by Theorem \ref{thm.++}.


\section{C-construction} The C-construction defined below is essential to study $\ZZ$-algebras.  We collect some properties of the C-construction which will be needed in this paper  (see \cite {Mbc}).

\subsection{C-construction}

\begin{definition} Let $\sC$ be a $k$-linear category.
\begin{enumerate}
\item{} ({\it B-construction}) For $\cO\in \sC$ and $s\in \Aut_k\sC$, we define a graded algebra $B(\sC, \cO, s):=\oplus_{i\in \ZZ}\Hom_{\sC}(\cO, s^i\cO)$.
\item{} ({\it C-construction}) For a sequence of objects $\{E_i\}_{i\in\ZZ}$ in $\sC$, we define a $\ZZ$-algebra $C(\sC, \{E_i\}_{i\in \ZZ}):=\oplus _{i, j\in \ZZ}\Hom_{\sC}(E_{-j}, E_{-i})$.
\end{enumerate}
\end{definition}

The following three lemmas are easy to prove, so we omit their proofs.

\begin{lemma} \label{lem.bc} \textnormal{\cite[Lemma 3.3]{Mbc}} Let $\sC$ be a $k$-linear abelian category.   For $\cO\in \sC$ and $s\in \Aut_k\sC$, $C(\sC, \{s^i\cO\}_{i\in \ZZ})\cong \overline {B(\sC, \cO, s)}$ as $\ZZ$-algebras so that
$$
\GrMod C(\sC, \{s^i\cO\}_{i\in \ZZ})\cong \GrMod B(\sC, \cO, s).
$$
\end{lemma}

\begin{lemma} \label{lem.czc} If $C$ is a $\ZZ$-algebra, then $C(r)\cong C(\GrMod C, \{P_{-i-r}\}_{i\in \ZZ})$ for every $r\in \ZZ$.
\end{lemma}

\begin{lemma} \label{lem.cec'e'} \textnormal{\cite[Lemma 2.5]{CN}}
Let $\{E_i\}_{i\in \ZZ}$ and $\{E'_i\}_{i\in \ZZ}$ be sequences in $k$-linear categories $\sC$ and $\sC'$, respectively.
\begin{enumerate}
\item{} A functor $F:\sC\to \sC'$ induces a $\ZZ$-algebra homomorphism
$$
C(\sC, \{E_i\}_{i\in \ZZ})\to C(\sC',  \{F(E_i)\}_{i\in \ZZ}).
$$
\item{} If there exists a fully faithful functor $F:\sC\to \sC'$ such that $F(E_i)\cong E'_i$ for every $i\in \ZZ$, then $C(\sC, \{E_i\}_{i\in \ZZ})\cong C(\sC', \{E'_i\}_{i\in \ZZ})$ as $\ZZ$-algebras.
\end{enumerate}
\end{lemma}

\begin{remark} \label{rem.end}
For a $\ZZ$-algebra $C$, the quotient functor $\pi : \GrMod C\to \Tails C$ induces a homomorphism
\begin{eqnarray*}
C & \cong & \uHom_C(C, C) \\
& = & C(\GrMod C, \{P_{-i}\}_{i\in \ZZ}) \\
& \to & C(\Tails C, \{\cP_{-i}\}_{i\in \ZZ}) \\
& = & \uHom_{\cC}(\cC, \cC)
\end{eqnarray*}
of $\ZZ$-algebras by Lemma \ref{lem.czc} and Lemma \ref{lem.cec'e'}, so $C(\Tails C, \{\cP_{-i}\}_{i\in \ZZ})$ has a bigraded $C$-bimodule structure.  On the other hand, $\uHom_{\cC}(\cC, \cC)$ has a bigraded $C$-bimodule structure as Definition \ref{def.mst}, and we can see that
$$
C(\Tails C, \{\cP_{-i}\}_{i\in \ZZ})=\uHom_{\cC}(\cC, \cC)
$$
as bigraded $C$-bimodules.
\end{remark}


\subsection{Ampleness}

The following notion of ampleness is a reindexed version of that introduced in \cite {P}.

\begin{definition}
We say that a sequence of objects $\{E_i\}_{i\in\ZZ}$ in an abelian category $\sC$ is {\it ample} if
\begin{enumerate}
\item[(A1)] for every $X\in \sC$ and every $m\in \ZZ$, there exists a surjection $\oplus _{j=1}^sE_{-i_j}\to X$ in $\sC$ for some $i_1, \dots, i_s\geq m$, and
\item[(A2)] for every surjection $X\to Y$ in $\sC$, there exists $m\in \ZZ$ such that $\Hom_{\sC}(E_{-i}, X)\to \Hom_{\sC}(E_{-i}, Y)$ is surjective for every $i \geq m$.
\end{enumerate}
\end{definition}

\begin{remark} \label{rem.ampZ}
For an abelian category $\sC$, $\cO\in \sC$ and $s\in \Aut \sC$, the pair $(\cO, s)$ is ample for $\sC$ in the sense of \cite {AZ} if and only if the sequence $\{s^i\cO\}_{i\in \ZZ}$ is ample for $\sC$ in the above sense.
\end{remark}

\begin{theorem} \label{thm.amp0}
Let $\sC$ be a $\Hom$-finite $k$-linear abelian category.  If $\{E_i\}_{i\in \ZZ}$ is an ample sequence for $\sC$ such that $\End_{\sC}(E_i)=k$ for every $i\in \ZZ$, then $C:=C(\sC, \{E_i\})_{\geq 0}$ is a right coherent connected $\ZZ$-algebra and there exists an equivalence functor $F:\sC\to \tails C$ such that $F(E_{-i})\cong \cP_i$ for every $i\in \ZZ$.
\end{theorem}

\begin{proof}
This follows from \cite[Theorem 2.3]{Na}, which is just a re-indexed version of \cite [Theorem 2.4]{P}.
\end{proof}

\subsection{$\chi$-condition}

\begin{definition} Let $C$ be a right coherent $\ZZ$-algebra.  We say that $C$  {\it satisfies $\chi_i$} if $\R^q\t(M)$ is right bounded for every $M\in \grmod C$ and every $q\leq i$.   We say that $C$ {\it satisfies $\chi$} if $C$ satisfies $\chi _i$ for every $i\in \NN$.
\end{definition}

\begin{remark} \label{rem.chi}
Let $C$ be a right $\Ext$-finite connected $\ZZ$-algebra.  By Lemma \ref{lem.tto}, $C$ satisfies $\chi_1$ if and only if, for every $M\in \grmod C$, there exists $m\in \ZZ$ such that the canonical map $M_{\geq m}\to (\omega \pi M)_{\geq m}$ is an isomorphism.
\end{remark}

The following lemma answers \cite[Remarks, page 70]{P}.

\begin{lemma} \label{lem.amp}
If $C$ is a right coherent connected $\ZZ$-algebra satisfying $\chi_1$, then $\tails C$ is $\Hom$-finite and $\{\cP_{-i}\}_{i\in \ZZ}$ is an ample sequence for $\tails C$.
\end{lemma}

\begin{proof}
 (A1):  For every $M\in \grmod C$, there exists a surjection  $F\to M$ in $\grmod C$ where $F\in \add\{P_j\}_{j\in \ZZ}$
which induces a surjection $\cF\to \cM$ in $\tails C$ where $\cF\in \add\{\cP_j\}_{j\in \ZZ}$,  so it is enough to show the condition (A1) for $\cP_i$ for  every $i\in \ZZ$.  Since $C$ is right $\Ext$-finite, for every $i\in \ZZ$, there exists an exact sequence $F\to P_i\to S_i\to 0$ in $\grmod C$ where $F\in \add\{P_j\}_{j>i}$ which induces a surjection $\cF\to \cP_i$ in $\tails C$ where $\cF\in \add\{\cP_j\}_{j>i}$.  Applying this argument finite number of times, for every $m\in \ZZ$, there exists a surjection $\cF\to \cP_i$ in $\tails C$ where $\cF\in \add\{\cP_j\}_{j\geq m}$.

(A2):
Every surjection $\pi \phi:\pi M\to \pi N$ in $\tails C$ where $M, N\in \grmod C$ is induced by a homomorphism $\phi:M\to N$ such that $\phi_{\geq m_1}:M_{\geq m_1}\to N_{\geq m_1}$ is a surjection in $\grmod C$ for some $m_1\in \ZZ$.  Since $C$ satisfies $\chi_1$, there exists $m_2, m_3\in \ZZ$ such that $M_{\geq m_2}\cong (\omega\pi M)_{\geq m_2}$ and $N_{\geq m_3}\cong (\omega\pi N)_{\geq m_3}$ by Remark \ref{rem.chi}.  By Lemma \ref{lem.Hom} (1), we have a commutative diagram
$$\begin{CD}
\Hom_{\cC}(\pi P_i, \pi M) @>\Hom_{\cC}(\pi P_i, \pi \phi)>> \Hom_{\cC}(\pi P_i, \pi N) \\
@V\cong VV @VV\cong V \\
\Hom_C(P_i, \omega\pi M) @>>> \Hom_C(P_i, \omega\pi N) \\
@V\cong VV @VV\cong V \\
(\omega\pi M)_i @>>>  (\omega\pi N)_i \\
@AAA @AAA \\
M_i
@>\phi_i>> N_i,
\end{CD}$$
so $\Hom_{\cC}(\pi P_i, \pi \phi):\Hom_{\cC}(\pi P_i, \pi M) \to \Hom_{\cC}(\pi P_i, \pi N)$ is surjective for every $i\geq \max \{m_1, m_2, m_3\}$.

We finally prove that $\tails C$ is $\Hom$-finite.  Let $\cM=\pi M, \cN=\pi N\in \tails C$ where $M, N\in \grmod C$.  By Remark \ref{rem.chi}, there exists $m\in \ZZ$ such that $N_{\geq m}\cong (\omega \pi N)_{\geq m}$.  By (A1), there exists a surjection $\oplus _{j=1}^s\cP_{i_j}\to \cM$ where $i_j\geq m$ for every $j=1, \dots, s$.
By Lemma \ref{lem.Hom} (1),
\begin{align*}
\Hom_{\cC}(\cM, \cN) & \subset \Hom_{\cC}(\oplus _{j=1}^s\pi P_{i_j}, \pi N) \cong  \oplus _{j=1}^s\Hom_C(P_{i_j},  \omega \pi N) \\
& =\oplus _{j=1}^s(\omega \pi N)_{i_j}=\oplus _{j=1}^sN_{i_j}.
\end{align*}
Since $C$ is locally finite by Lemma \ref{lem.lofi}, $N$ is locally finite, so $\Hom_{\cC}(\pi M, \pi N)$ is finite dimensional.
\end{proof}

\begin{lemma} \label{lem.aschi}
If $C$ is a right coherent ASF-regular $\ZZ$-algebra of dimension $d\geq 1$, then $C$ satisfies $\chi$, so  $\tails C$ is $\Hom$-finite and $\{\cP_{-i}\}_{i\in \ZZ}$ is an ample sequence for $\tails C$.
\end{lemma}

\begin{proof} Since $C$ is right coherent, for every $M\in \grmod C$, there exist $F\in \add \{P_{i}\}_{i \in \mathbb{Z}}$ and $L\in \grmod C$ such that $0\to L\to F\to M\to 0$ is an exact sequence, which induces an exact sequence $\R^q\t(F)\to \R^q(M)\to \R^{q+1}\t(L)$ in $\GrMod C$ for every $q\in \ZZ$.  Since $C$ is ASF-regular, $\R^q\t(P_i)$ is right bounded for every $i\in \ZZ$ and $q\in \ZZ$, so $\R^q\t(F)$ is right bounded for every $q\in \ZZ$.  Since $\pd(M)\leq d$ by \cite[Proposition 4.11]{MN},
$\R^q\t(M)=\R^q\t(L)=0$ are right bounded for every $q>d$.  In particular, $\R^{d+1}\t(L)=0$, so $\R^d\t(M)$ is right bounded.  By the same argument, $\R^d\t(L)$ is right bounded, so $\R^{d-1}\t(M)$ is right bounded.  By induction, $\R^q\t(M)$ is  right bounded for every $q\in \ZZ$, so $C$ satisfies $\chi$.  By Lemma \ref{lem.amp}, $\tails C$ is $\Hom$-finite and $\{\cP_{-i}\}_{i\in \ZZ}$ is an ample sequence for $\tails C$.
\end{proof}

\subsection{Quasi-Veronese $\ZZ$-algebras}

\begin{definition} Let $I_j:=\{i\in \ZZ\mid jr\leq i\leq (j+1)r-1\}$ for each $j\in \ZZ$.  For a $\ZZ$-algebra $C$ and $r\in \NN^+$,
the {\it $r$-th quasi-Veronese $\ZZ$-algebra $C^{[r]}$} of $C$ is defined by $(C^{[r]})_{st}:=\oplus _{-i\in I_{-s}, -j\in I_{-t}}C_{ij}$
\end{definition}

Note that if $C=C(\sC, \{E_i\}_{i\in \ZZ})$, then $C^{[r]}=C(\sC, \{\oplus _{i\in I_j}E_i\}_{j\in \ZZ})$.

\begin{lemma} \label{lem.qv}
For a $\ZZ$-algebra $C$ and $r\in \NN^+$, $\GrMod C\cong \GrMod C^{[r]}$.
\end{lemma}

\begin{proof}  Since $C^{[r]}\cong C$ as algebras (not as $\ZZ$-algebras unless $r=1$) by a variant of \cite[Corollary 3.6]{Mbc}, the result holds by Lemma \ref{lem.mgm2}.
\end{proof}

\begin{lemma} \label{lem.aqv}
A sequence of objects $\{E_i\}_{i\in \ZZ}$ in an abelian category $\sC$ is ample if and only if $\{\oplus _{i\in I_j}E_i\}_{j\in \ZZ}$ is ample.
\end{lemma}

\begin{proof}
(A1) Let $M\in \sC$.  If $\{E_i\}_{i\in \ZZ}$ is ample, then, for every $m\in \ZZ$, there exists a surjection $\oplus _{q=1}^{s}E_{-i_q}\to M$ in $\sC$ for some $i_1, \dots, i_{s}\geq rm-(r-1)$.  Since, for $1 \leq q \leq s$, $i_{q}\geq rm-(r-1)$, we have $-i_{q} \leq -rm+(r-1)$, and thus there exists $\ell_q \in \mathbb{Z}$ such that $\ell_q \geq m$ and
$$
-r \ell_q \leq -i_{q} \leq -r \ell_{q}+(r-1).
$$
Thus, $-i_{q} \in I_{-\ell_q}$ and therefore there exists a surjection
$$
\bigoplus _{q=1}^{s}\bigoplus_{i \in I_{-\ell_{q}}}E_{i}\to M.
$$

Conversely, if $\{\oplus _{i\in I_j}E_i\}_{j\in \ZZ}$ is ample, then, for every $m\in \ZZ$, there exists a surjection $\oplus _{q=1}^{s}(\oplus _{i\in I_{-j_q}}E_{i})\to M$ in $\sC$ for some $j_1, \dots, j_{s}\geq (m+r-1)/r$.  Thus, for $1 \leq q \leq s$, $rj_{q} \geq m+(r-1)$ and thus $-rj_{q}+(r-1) \leq -m$.  It follows that if $i \in I_{-j_{q}}$, then $-i \geq m$, and so $\{E_i\}_{i\in \ZZ}$ is ample.

(A2) Let $M\to N$ be a surjection in $\sC$.  If $\{E_i\}_{i\in \ZZ}$ is ample, then there exists $m\in \ZZ$ such that $\Hom_{\sC}(E_{-i}, M)\to \Hom_{\sC}(E_{-i}, N)$ is surjective for every $i\geq m$.  Now suppose $\ell \geq (m+r-1)/r$.  Then $r \ell - (r-1) \geq m$ and thus if $w \in I_{-\ell}$ then $m \leq -w$.  It follows that $\Hom_{\sC}(\oplus _{i\in I_{-\ell}}E_i, M)\to \Hom_{\sC}(\oplus _{i\in I_{-\ell}}E_i, N)$ is surjective for every $\ell \geq (m+r-1)/r$.

Conversely, if $\{\oplus _{i\in I_j}E_i\}_{j\in \ZZ}$ is ample, then there exists $m\in \ZZ$ such that $\Hom_{\sC}(\oplus _{i\in I_{-j}}E_i, M)\to \Hom_{\sC}(\oplus _{i\in I_{-j}}E_i, N)$ is surjective for every $j \geq m$.

Now suppose $-i \geq mr-(r-1)$.  Then $i \leq -mr+(r-1)$ and thus $i \in I_{-j}$ for some $j \geq m$ so that $\Hom_{\sC}(E_{i}, M)\to \Hom_{\sC}(E_{i}, N)$ is surjective.
\end{proof}


\subsection{Helices}

In this subsection, we recall some notions from \cite{MU1}.

\begin{definition} Let $\sC$ be a $\Hom$-finite $k$-linear category.
A {\it Serre functor for $\sC$} is a $k$-linear autoequivalence $S\in \Aut _k\sC$ such that there exists a bifunctorial isomorphism
$$F_{X, Y}:\Hom_{\sC}(X, Y)\to D\Hom_{\sC}(Y, S(X))$$ for $X, Y\in \sC$.
\end{definition}

\begin{remark} \label{rem.Se2}
We explain the functoriality of a Serre functor $S$ in $Y$ in the above definition.   (See \cite[Remark 3.2]{MU1} for the functoriality in $X$.) Define functors $G=\Hom_{\sC}(X, -)$ and $H=D\Hom_{\sC}(-, S(X))=\Hom_k(\Hom_{\sC}(-, S(X)), k)$.  Fix $\b\in \Hom_{\sC}(Y, Y')$.  Then
$$G(\b):\Hom_{\sC}(X, Y)\to \Hom_{\sC}(X, Y')$$
is given by $(G(\b))(\a)=\b\circ \a$.  On the other hand,
$$H(\b):\Hom_k(\Hom_{\sC}(Y, S(X)), k)\to \Hom_k(\Hom_{\sC}(Y', S(X)), k)$$
is given by $((H(\b))(\phi))(\c)=\phi(\c\circ \b)$ for $\c\in \Hom_{\sC}(Y', S(X))$.
By functoriality, we have a commutative diagram
$$\begin{CD}
\Hom_{\sC}(X, Y) @>{G(\b)}>> \Hom_{\sC}(X, Y') \\
@VF_{X, Y}VV @VVF_{X, Y'}V \\
D\Hom_{\sC}(Y, S(X)) @>{H(\b)}>> D\Hom_{\sC}(Y', S(X)),
\end{CD}$$
so, for $\a\in \Hom_{\sC}(X, Y)$ and $\c\in \Hom_{\sC}(Y', S'(X))$, we have
$$F_{X, Y'}(\b\circ \a)(\c)=(F_{X, Y'}(G(\b)(\a)))(\c)=(H(\b)(F_{X, Y}(\a)))(\c)=F_{X, Y}(\a)(\c\circ \b).$$
\end{remark}

\begin{definition} Let $\sC$ be an abelian category.  A {\it bimodule $\cL$ over $\sC$} is an adjoint pair of functors from $\sC$ to itself with the suggestive notation $\cL=(-\otimes_{\sC}\cL, \cHom_{\sC}(\cL, -))$.  A bimodule $\cL$ over $\sC$ is {\it invertible} if $-\otimes _{\sC}\cL$ is an autoequivalence of $\sC$.  In this case, the inverse bimodule of $\cL$ is defined by $\cL^{-1}=(-\otimes _{\sC}\cL^{-1}, \cHom_{\sC}(\cL^{-1}, -))=(\cHom_{\sC}(\cL, -), -\otimes _{\sC}\cL)$.
\end{definition}


\begin{definition} \label{dfn.cf}
Let $\sC$ be a $\Hom$-finite $k$-linear abelian category.
A {\it canonical bimodule for $\sC$} is
an invertible bimodule $\omega _{\sC}$ over $\sC$ such that, for some $n\in \ZZ$,
the autoequivalence $-\lotimes _{\sC}\omega _{\sC}[n]$ of $\sD^b(\sC)$ induced by $-\otimes _{\sC}\omega _{\sC}$ is a Serre functor for $\sD^b(\sC)$.
\end{definition}

\begin{definition} \label{def.BV}
Let $\sT$ be a triangulated category. For a set of objects
$\cE$ in $\sT$, we denote by
$\<\cE\>$ the smallest full triangulated subcategory of $\sT$ containing
$\cE$ closed under isomorphisms and direct summands.
We say that $S$ classically generates $\sT$ if $\<\cE\>=\sT$.
\end{definition}

\begin{definition} \label{dfn.res}
Let $\sT$ be a $k$-linear triangulated category, and $\cE=\{E_0, \dots, E_{\ell-1}\}$ a sequence of objects in $\sT$.
\begin{enumerate}
\item{}
$\cE$ is called {\it exceptional} if
\begin{enumerate}
\item[(E1)] $\End_{\sT}(E_i)=k$
for every $i=0, \dots, \ell-1$,
\item[(E2)] $\Hom_{\sT}(E_i, E_i[q])=0$ for every $q\neq 0$ and every $i=0, \dots, \ell-1$, and
\item[(E3)] $\Hom_{\sT}(E_i, E_j[q])=0$ for every $q$ and every $0\leq j<i\leq \ell-1$.
\end{enumerate}
\item{}
$\cE$ is called {\it full} if
$\<\cE\>=\sT$.
\end{enumerate}
\end{definition}

\begin{definition} \label{dfn.rh}
Let $\sC$ be a $k$-linear abelian category having the canonical bimodule $\omega _{\sC}$, and $\cE=\{E_i\}_{i\in \ZZ}$ a sequence of objects in $\sD^b(\sC)$.
\begin{enumerate}
\item{}
$\cE$ is called a
{\it helix of period $\ell$} if
\begin{enumerate}
\item[(H1)] $E_{i+\ell}\cong E_i\lotimes _{\sC}\omega _{\sC}^{-1}$ for every $i$,
\item[(H2)] $\Ext^q_{\sC}(E_i, E_i)\cong \begin{cases} k & \textnormal { if } q=0, \\
0 & \textnormal { if } q\neq 0\end{cases}$ for every $i$,  and
\item[(H3)] $\Ext^q_{\sC}(E_i, E_j)=0$ for every $q$ and every $0<i-j<\ell$
(or equivalentry every $j<i<j+\ell$).
\end{enumerate}
\item{} A
helix $\cE$ of period $\ell$ is called {\it geometric} if
$\Ext^q_{\sC}(E_i, E_j)=0$
for every $q\neq 0$ and every $i\leq j$.
\item{} A
helix $\cE$ of period $\ell$ is called {\it full} if $\<E_{i}, \dots, E_{i+\ell-1}\>=\sD^b(\sC)$ for every $i$.
\end{enumerate}
\end{definition}

\begin{remark} Let $\sC$ be a $k$-linear abelian category having the canonical bimodule $\omega _{\sC}$, and $\cE:=\{E_i\}_{i\in \ZZ}$ a
(full) helix of period $\ell$ for $\sD^b(\sC)$.
\begin{enumerate}
\item{} $\{E_{i}, \dots, E_{i+\ell-1}\}$ is a
(full) exceptional sequence for $\sD^b(\sC)$ for every $i$.
\item{}  If $\gldim \sC=n$, then the Serre functor for $\sD^b(\sC)$ is given by $-\lotimes _{\sC}\omega _{\sC}[n]$ by \cite [Remark 3.5]{MU1}, so
$$\Hom_{\sC}(E_i, E_j)\cong D\Hom_{\sC}(E_j, E_i\lotimes _{\sC}\omega_{\sC}[n])\cong D\Ext^{n}_{\sC}(E_j, E_{i-\ell})$$
for every $i, j\in \ZZ$.
\item{} The above definition of a helix is not exactly the same as the one defined by mutation in \cite{BP} (see \cite[Remark 3.17]{MU1}).
\end{enumerate}
\end{remark}


\begin{lemma} \label{lem.olp}
If $\sC$ is a $k$-linear abelian category having the canonical bimodule $\omega _{\sC}$, and $\{E_i\}_{i\in \ZZ}$ is a helix of period $\ell$ for $\sD^b(\sC)$, then $C=C(\sC, \{E_i\}_{i\in \ZZ})$ is $\ell$-periodic.
\end{lemma}

\begin{proof}
Since $\omega_{\sC}^{-1}:\sC\to \sC$ is an equivalence functor such that $E_i\otimes \omega _{\sC}^{-1}\cong E_{i+\ell}$ for every $i\in \ZZ$, $C(\ell)=C(\sC, \{E_{i+\ell}\}_{i\in \ZZ})\cong C(\sC, \{E_i\}_{i\in \ZZ})=C$ as $\ZZ$-algebras by Lemma \ref{lem.cec'e'}, so $C$ is $\ell$-periodic.
\end{proof}


\section{A Categorical Characterization of Quantum Projective $\ZZ$-spaces}

The main theorem of \cite{MU1} is a categorical characterization of quantum projective spaces, which are defined to be the noncommutative projective schemes associated to AS-regular algebras.  In this last section, we will give a $\ZZ$-algebra version of this result.  As a biproduct, we give a family of non-trivial examples of AS-regular $\ZZ$-algebras, which are constructed from noncommutative quadric hypersurfaces.

\subsection{A Categorical Characterization}


\begin{definition}
Let $\sC$ be a $\Hom$-finite $k$-linear abelian category.  We say that $\sC$ satisfies {\it (GH) of period $\ell$} if
\begin{enumerate}
\item[(GH1)] $\sC$ has a canonical bimodule $\omega _{\sC}$, and
\item[(GH2)] $\sC$ has an ample sequence $\{E_i\}_{i\in \ZZ}$ which is a full geometric helix  of period $\ell$ for $\sD^b(\sC)$.
\end{enumerate}
\end{definition}

Let $\sC$ be a $\Hom$-finite $k$-linear abelian category satisfying (GH) of period  $\ell$.  We write
\begin{itemize}
\item{} $T:=E_0\oplus \cdots \oplus E_{\ell-1}
\in \sC$,
\item{} $R:=\End_{\sC}(T)$, and
\item{} $\Pi R:=B(\sC, T, -\otimes _{\sC}\omega _{\sC}^{-1})_{\geq 0}$.
\end{itemize}

\begin{lemma} \label{lem.pir}
If $\sC$ is a $\Hom$-finite $k$-linear abelian category satisfying (GH) of period  $\ell$, then $\gldim \Pi R=\gldim \sC+1<\infty$.
\end{lemma}

\begin{proof}  Since $\{T\otimes _{\sC}(\omega _{\sC}^{-1})^{\otimes j}\}_{j\in \ZZ}=\{\oplus _{i\in I_j}E_i\}_{j\in \ZZ}$ where
$$
I_j=\{i\in \ZZ\mid j\ell\leq i\leq (j+1)\ell-1\}
$$
is ample for $\sC$ by Lemma \ref{lem.aqv},
$(T, -\otimes_{\sC}\omega _{\sC}^{-1})$ is ample for $\sC$ in the sense of \cite {AZ} by Remark \ref{rem.ampZ}.
By \cite [Lemma 3.18, Lemma 3.19, Theorem 3.11 (2)]{MU1},
$\gldim \Pi R=\gldim \sC+1<\infty$.
\end{proof}

\begin{remark} \label{rem.reti}  In the above lemma,  $\{T\otimes _{\sC}(\omega _{\sC}^{-1})^{\otimes i}\}_{i\in \ZZ}$ is a ``full geometric relative helix of period 1" for $\sD^b(\sC)$ in the sense of \cite [Definition 3.14]{MU1} by \cite[Lemma 3.18]{MU1}, so  $T$ is  a ``regular tilting object'' for $\sD^b(\sC)$ in the sense of \cite [Definition 3.9]{MU1} by \cite[Lemma 3.19]{MU1}.
In the literature, $R=\End_{\sC}(T)$ is called a ``Fano algebra'', and $\Pi R:=B(\sC, T, -\otimes _{\sC}\omega _{\sC}^{-1})_{\geq 0}$ is called the ``preprojective algebra" of $R$ (cf. \cite{MM}).
\end{remark}




The following theorem is a generalization of \cite[Theorem 4.2]{BP} and one direction of  \cite[Theorem 4.1]{MU1}.
Our definition of a helix is not the same as the one defined in \cite{BP},  so that the Koszul condition plays no role in the theorem below (see \cite[Remark 3.17, Corollary 4.2]{MU1}).

\begin{theorem} \label{thm.cs1}

If $\sC$ is a $\Hom$-finite $k$-linear abelian category satisfying (GH) of period  $\ell$, then $C:=C(\sC, \{E_i\}_{i\in \ZZ})_{\geq 0}$ is a
right coherent ASF$^{++}$-regular $\ZZ$-algebra of dimension $\gldim \sC+1$ and of Gorenstein parameter $\ell$ such that $\sC\cong \tails C$.
\end{theorem}

\begin{proof}
Suppose that $\sC$ satisfies (GH) of period $\ell$ with $n=\gldim \sC<\infty$.  Since $\sC$ is $\Hom$-finite, $C$ is locally finite.  Since $C_{ii}=\Hom_{\sC}(E_i, E_i)=k$ by (H2), $C$ is a connected $\ZZ$-algebra.
 By Theorem \ref{thm.amp0}, $C$ is right coherent and $\sC\cong \tails C$, so it is enough to show that $C$ is an $\ell$-periodic ASF$^+$-regular $\ZZ$-algebra of dimension $n+1$ and of Gorenstein parameter $\ell$ by Theorem \ref{thm.++}.
Let $T:=E_0\oplus \cdots \oplus E_{\ell-1}\in \sC$ and $R:=\End_{\sC}(T)$.

\smallskip

First suppose that $n=0$.  Since $\{E_0, \dots, E_{\ell-1}\}$ is an exceptional sequence, $R=\End_{\sC}(\oplus _{i=0}^{\ell-1}E_i)=\begin{pmatrix} k & \Hom_{\sC}(E_0, E_1) & \cdots & \Hom_{\sC}(E_0, E_{\ell-1}) \\
0 & k & \cdots & \Hom_{\sC}(E_1, E_{\ell-1}) \\
\vdots & \vdots & \ddots & \vdots \\
0 & 0 & \cdots & k\end{pmatrix}$ is an upper triangular matrix algebra over $k$.  By Remark \ref{rem.reti} and \cite[Theorem 3.10]{MU1}, $R$ is semisimple, so $R\cong k^{\ell}$, hence $\Hom_{\sC}(E_i, E_j)=0$ for any $0\leq i\neq j\leq \ell-1$.   Since $\cD^b(\sC)\cong \cD^b(\mod R)$, we have $\omega _{\sC}=\id _{\sC}$, so $E_{i+s\ell}\cong E_i$ for every $s\in \ZZ$.  It follows that, for every $i, j\in \ZZ$, there exists $0\leq i', j'\leq \ell-1$ such that $E_i\cong E_{i'}, E_j\cong E_{j'}$,
so
\begin{align*}
C(\sC, \{E_i\}_{i\in \ZZ})_{ij} & \cong \Hom_{\sC}(E_{-j}, E_{-i})\cong \Hom_{\sC}(E_{-j'}, E_{-i'}) \\
& \cong \begin{cases} k & \textnormal { if $i'=j'$, or equivalently if $\ell \;\;  | \; j-i$} \\ 0 & \textnormal { if $i'\neq j'$, or equivalently if $\ell \not| \; j-i$,} \end{cases} \\
& \cong k[x, x^{-1}]_{j-i}\cong \overline {k[x, x^{-1}]}_{ij}
\end{align*}
where $\deg x=\ell$.  Since we have a canonical commutative diagram
$$\begin{CD}
C(\sC, \{E_i\}_{i\in \ZZ})_{i,i+s\ell} \otimes C(\sC, \{E_i\}_{i\in \ZZ})_{i+s\ell, i+(s+t)\ell} @>>> C(\sC, \{E_i\}_{i\in \ZZ})_{i,i+(s+t)\ell} \\
@V\cong VV @VV\cong V \\
\Hom_{\sC}(E_{-(i+s\ell)}, E_{-i}) \otimes \Hom_{\sC}(E_{-(i+(s+t)\ell)}, E_{-(i+s\ell)})  @>>> \Hom_{\sC}(E_{-(i+(s+t)\ell)}, E_{-i}) \\
@V\cong VV @VV\cong V \\
\Hom_{\sC}(E_{-i}, E_{-i})\otimes \Hom_{\sC}(E_{-i}, E_{-i}) @>>> \Hom_{\sC}(E_{-i}, E_{-i}) \\
@V\cong VV @VV\cong V \\
k \otimes k @>>> k \\
@V\cong VV @VV\cong V \\
kx^s \otimes kx^t @>>> kx^{s+t} \\
@V\cong VV @VV\cong V \\
\overline {k[x, x^{-1}]}_{i, i+s\ell} \otimes \overline {k[x, x^{-1}]}_{i+s\ell, i+(s+t)\ell} @>>> \overline {k[x, x^{-1}]}_{i, i+(s+t)\ell}
\end{CD}$$
$C:=C(\sC, \{E_i\}_{i\in \ZZ})_{\geq 0}\cong \overline {k[x, x^{-1}]}_{\geq 0}\cong \overline {k[x]}$ as $\ZZ$-algebras.  Since $k[x]$ is an ASF-regular algebra of dimension 1 and of Gorenstein parameter $\ell$, $C$ is an $\ell$-periodic (in fact 1-periodic) ASF$^+$-regular $\ZZ$-algebra of dimension 1 and of Gorenstein parameter $\ell$ by Theorem \ref{thm.asfo}.

\smallskip

Now suppose that $n\geq 1$.  If $j<i<j+\ell$, then $\Hom_{\sC}(E_i, E_j)=0$ by (H3).
If $j+\ell\leq i$, then $\Hom_{\sC}(E_i, E_j)\cong D\Ext^n_{\sC}(E_j, E_{i-\ell})=0$ since $\{E_i\}$ is geometric and $n\geq 1$.  It follows that $C=C(\sC, \{E_i\}_{i\in \ZZ})$, so $C$ is $\ell$-periodic by Lemma \ref{lem.olp}.

\smallskip

\noindent (ASF1):  Since $E_i\lotimes _{\sC}\omega _{\sC}^{-1}\cong E_{i+\ell}$ for every $i\in \ZZ$ by (H1), and
$\Hom_{\sC}(E_i, E_j)=0$ for every $j<i$ as shown above,
$$\Pi R := B(\sC, T, -\otimes _{\sC}\omega _{\sC}^{-1})_{\geq 0}=B(\sC, T, -\otimes _{\sC}\omega _{\sC}^{-1}).$$  Since $T\otimes _{\sC}(\omega _{\sC}^{-1})^{\otimes j}=\oplus _{i\in I_j}E_i$ where $I_j=\{i\in \ZZ\mid j\ell\leq i\leq (j+1)\ell-1\}$,
\begin{eqnarray*}
\overline {\Pi R}
& = & \overline {B(\sC, T, -\otimes _{\sC}\omega _{\sC}^{-1})} \\
& \cong & C(\sD^b(\sC), \{
T\lotimes _{\sC}\omega _{\sC}^{-j}\}_{j\in \ZZ}) \\
& \cong & C(\sD^b(\sC), \{\oplus _{i\in I_j} E_i\}_{j\in \ZZ}\}) \\
& \cong & C(\sC, \{E_i\}_{i\in \ZZ})^{[\ell]} \\
& = & C^{[\ell]},
\end{eqnarray*}
so $\GrMod C\cong \GrMod C^{[\ell]}\cong \GrMod \overline {\Pi R}\cong \GrMod \Pi R$ by Lemma \ref{lem.qv} and Lemma \ref{lem.ovwi}.
Since $\{S_i\}_{i\in \ZZ}$ is the set of complete representatives of isomorphism classes of simple objects in $\GrMod C$,
$$\operatorname {sgldim} C=\operatorname {sgldim}\Pi R=\gldim \Pi R=n+1$$
by \cite [Proposition 2.7]{MM} and Lemma \ref{lem.pir}.


\smallskip

\noindent (ASF2), (ASF2$^+$):
We will check the equivalent conditions given in Proposition \ref{prop.asf2'}.   Since $\{E_i\}_{i\in \ZZ}$ is ample for $\sC$, there exists an equivalence functor $\sC\to \tails C$
 sending $E_{-i}$ to $\cP_i$ for every $i\in \ZZ$ by Theorem \ref{thm.amp0}.   Since $\uHom_{\cC}(\cC, \cC)=C(\tails C, \{\cP_{-i}\}_{i\in \ZZ})\cong C(\sC, \{E_i\}_{i\in \ZZ})=: C$ as bigraded $C$-bimodules by Lemma \ref{lem.cec'e'} and Remark \ref{rem.end}, $e_i\uHom_{\cC}(\cC, \cC)\cong e_iC=P_i$ as graded right $C$-modules for every $i\in \ZZ$ and $\uHom_{\cC}(\cC, \cC)e_j\cong Ce_j=Q_j$ as graded left $C$-modules for every $j\in \ZZ$.

If $i\leq j$, then $\Ext^q_{\cC}(\cP_j, \cP_i)\cong \Ext^q_{\sC}(E_{-j}, E_{-i})=0$ for all $q\neq 0$ since $\{E_i\}_{i\in \ZZ}$ is geometric.
If $i-\ell< j < i$, then $\Ext^q_{\cC}(\cP_j, \cP_i)\cong \Ext^q_{\sC}(E_{-j}, E_{-i})=0$ for all $q$ by (H3).
If $j\leq i-\ell$, then
$$
\Ext^q_{\cC}(\cP_j, \cP_i)\cong \Ext^q_{\sC}(E_{-j}, E_{-i})\cong D\Ext^{n-q}_{\sC}(E_{-i}, E_{-j-\ell})=0
$$
for all $q\neq n$ since $\{E_i\}_{i\in \ZZ}$ is geometric. It follows that
$$
\uExt^q_{\cC}(\cC, \cC)=\oplus _{i, j\in \ZZ}\Ext_{\cC}^q(\cP_j, \cP_i)=0
$$
unless $q=0, n$.

We will prove only (ASF2$^+$).   The proof for (ASF2) is similar (see the proof of \cite[Theorem 4.1]{MU1}).
Let $F_{ij}$ be the isomorphism
\begin{eqnarray*}
\Ext_{\sC}^{d-1}(E_{-j}, E_{-i}) & \cong & \Hom_{\sC}(E_{-j}[1-d], E_{-i}) \\
& \rightarrow & D\Hom_{\sC}(E_{-i}, S(E_{-j}[1-d])) \\
& \cong & D\Hom_{\sC}(E_{-i}, E_{-j}\lotimes _{\sC}\omega_{\sC}) \\
& \cong & D\Hom_{\sC}(E_{-i}, E_{-j-\ell})
\end{eqnarray*}
induced by the Serre functor, and consider the diagram
$$\begin{CD}
C_{ij}\times \Ext_{\cC}^{d-1}(\cC, \cC)_{jk} @>>> \Ext_{\cC}^{d-1}(\cC, \cC)_{ik} \\
@VVV @VVV \\
\Hom_{\sC}(E_{-j}, E_{-i})\times \Hom_{\sC}(E_{-k}[1-d], E_{-j}) @>\Phi>> \Hom_{\cC}(E_{-k}[1-d], E_{-i}) \\
@V\id\times F_{jk} VV @VV{F_{ik}}V \\
\Hom_{\sC}(E_{-j}, E_{-i})\times D\Hom_{\sC}(E_{-j}, E_{-k-\ell}) @>\Psi>> D\Hom_{\cC}(E_{-i}, E_{-k-\ell}) \\
@VVV @VVV \\
C_{ij}\times D(C(\ell,0))_{jk} @>>> D(C(\ell,0))_{ik}.
\end{CD}$$
The top square commutes since this is the way to define the graded left $C$-module structure for $\uExt_{\cC}^{d-1}(\cC, \cC)$ as in Definition \ref{def.mst} (2).  The middle square commutes by the bifunctoriality of the Serre functor as follows: For $(\b, \a)\in \Hom_{\sC}(E_{-j}, E_{-i})\times \Hom_{\sC}(E_{-k}[1-d], E_{-j})$, we have $\Phi(\b, \a)=\b\circ \a$.  Moreover, for $(\b, \phi)\in  \Hom_{\sC}(E_{-j}, E_{-i})\times D\Hom_{\sC}(E_{-j}, E_{-k-\ell})$, we have $\Psi (\b, \phi)(\c)=\phi (\c\circ \b)$ for every $\c\in \Hom_{\cC}(E_{-i}, E_{-k-\ell})$.  By Remark \ref{rem.Se2},
\begin{align*}
F_{ik}(\Phi(\b, \a))(\c) & =F_{ik}(\b\circ \a)(\c)=F_{jk}(\a)(\c\circ \b) \\
& =\Psi(\b, F_{jk}(\a))(\c)=\Psi((\id\times F_{jk})(\b, \a))(\c).
\end{align*}
The bottom square commutes since this is the way to define the graded left $C$-module structure for $D(C(\ell,0))$, so $\uExt_{\cC}^{d-1}(\cC, \cC)\cong D(C(\ell,0))$ as graded left $C$-modules, hence (ASF2$^+$) holds.
\end{proof}

\begin{theorem} \label{thm.npl}

Let $\sC$ be a $k$-linear abelian category.  If $\sC$ satisfies (GH), then the following are equivalent:
\begin{enumerate}
\item{} $C:=C(\sC, \{E_i\}_{i\in \ZZ})_{\geq 0}$ is right noetherian.
\item{} $\sC$ is a noetherian category.
\item{} $E_i\in \sC$ is a noetherian object for every $i\in \ZZ$.
\end{enumerate}
\end{theorem}

\begin{proof}
(1) $\Rightarrow$ (2): If $C$ is right noetherian, then
$\grmod C$ is a noetherian category by Lemma \ref{lem.noe}.
It follows that, for every $\cM\in \tails C$, there exists a noetherian object $M\in \grmod C$ such that $\cM\cong \pi M\in \tails C$. Since $C$ is a right coherent connected $\ZZ$-algebra by Lemma \ref{lem.cone},  
$\Tors C$ is a localizing subcategory of $\GrMod C$ by Lemma \ref{lem.loc}, so $\cM\cong \pi M\in \Tails C$ is a noetherian object by \cite[Lemma 5.8.3]{Po}, hence $\tails C$ is a noetherian category.

(2) $\Rightarrow$ (3): If $\sC$ is a noetherian category, then $E_i\in \sC$ is a noetherian object for every $i\in \ZZ$ by definition.

(3) $\Rightarrow$ (1): If $E_i\in \sC$ is a noetherian object for every $i\in \ZZ$, then $T:=E_1\oplus \cdots \oplus E_{\ell-1}\in\sC$ is a noetherian object, so $\Pi R$ is right noetherian by \cite [Theorem 4.1]{MU1}.  It follows that $\grmod C\cong \grmod \Pi R$ is a noetherian category, so $C$ is right noetherian by Lemma \ref{lem.noe}.
\end{proof}

We will next show  the converse to Theorem \ref{thm.cs1} to complete our categorical characterization (a generalization of the other direction of \cite[Theorem 4.1]{MU1}).

\begin{theorem} \label{thm.LD3}
Let $C$ be a right coherent ASF$^{++}$-regular $\ZZ$-algebra.  If $X, Y \in \sD^b(\grmod C)$, then we have a bifunctorial isomorphism
$$\R\uHom_C(I_{0}(X), I_{0}(Y)\ulotimes_C D\R Q(C))
\cong D\R\uHom_C(I_{0}(Y), \R Q(I_{0}(X))).$$
\end{theorem}

\begin{proof}
By Lemma \ref{lem.sgl2}, $I_{0}(X)$ and $I_{0}(Y)$ are quasi-isomorphic to perfect complexes.  Furthermore, since $X$ is bounded and $D\R Q(C)$ is bounded by Lemma \ref{lem.bound},
$D\R Q(I_{0}(X))\in \sD^b (\Bimod (C-K))$ by Theorem \ref{thm.ld0}. We next note that by Lemma \ref{lem.tto}, there is a triangle
$$
(D\R Q(X))_{j} \rightarrow (DX)_{j} \rightarrow (D \R \t (X))_{j}
$$
in $\sD(\Mod k)$.  Since, by Theorem \ref{thm.ld1}, $D \R \t (X)$ has locally finite homology, so does $D \R Q(X)$.  It follows that $D\R Q(I_{0}(X))$ has locally finite homology.  Therefore, by Proposition \ref{prop.identity}, Theorem \ref{thm.ld0} and \cite[Corollary 6.2]{MN}, we have bifunctorial isomorphisms
\begin{align*}
\R\uHom_C(I_{0}(X), I_{0}(Y) \ulotimes_C D\R Q(C))
& \cong I_{0}(Y) \ulotimes _C\R\uHom_C(I_{0}(X),  D\R Q(C)) \\
& \cong I_{0}(Y) \ulotimes_CD\R Q(I_{0}(X)) \\
& \cong D\R\uHom_C(I_{0}(Y), \R Q(I_{0}(X))),
\end{align*}
where in the last isomorphism, in order to apply \cite[Corollary 6.2]{MN}, we use \cite[Proposition 3.1(1)]{Ve}, which holds since the last two expressions have locally finite homology by \cite[Proposition 7.3(i)]{hartshorne} and the fact that $I_{0}(Y)$ is quasi-isomorphic to a perfect complex and $\R Q(I_{0}(X)))$ has locally finite homology. 
\end{proof}


\begin{lemma} \label{lem.lln}
If $C$ is a right coherent ASF$^{++}$-regular $\ZZ$-algebra and $X, Y\in \sD^b(\grmod C)$, then
\begin{enumerate}
\item{} if $Y_{\geq n}=0$ for some $n \in \ZZ$, then $\RuHom_C(X_{\geq n}, Y)=0$, and

\item{} $\RuHom_C(X_{\geq n}, \R Q(Y))\cong \R\uHom_C(X_{\geq n},  Y)$ for some $n\in \ZZ$.
\end{enumerate}
\end{lemma}

\begin{proof}
(1) Since $C$ is right coherent, $X/X_{\geq n}\in \sD^b(\grmod C)$, so $X_{\geq n}\in \sD^b(\grmod C)$.  If $F$ is the minimal free resolution of $X_{\geq n}$, then $F_{<n}=0$, so
$$
\RuHom_C(X_{\geq n},Y)=\uHom_C(F,Y)=0,
$$
hence the result.

(2)
For $X, Y\in \sD^b(\grmod C)$ and $n \in \mathbb{Z}$, we have a triangle $\R\t(Y)\to Y\to \R Q(Y)$ in $\sD(\GrMod C)$ by Lemma \ref{lem.tto}, which induces a triangle
$$\RuHom_C(X_{\geq n},\R\t(Y)) \to \RuHom_C(X_{\geq n}, Y)\to \RuHom_C(X_{\geq n},\R Q(Y)).$$

Since $C$ is right coherent and $\operatorname{sgldim}C<\infty$, $Y\in \sD^b(\grmod C)$ has a finitely generated free resolution $G$ of finite length by Lemma \ref{lem.sgl2}.  Since $C$ is ASF$^{++}$-regular, it is ASF-regular and so $\R\t(P_j)\cong (D(Q_{j-\ell}))[-d]$ is a complex whose terms are right  bounded for every $j\in \ZZ$, so $\R\t(G)_{\geq n}=0$ for some $n\in \ZZ$, and so the result now follows from (1).
\end{proof}

The proof of the following result was inspired by the proof of \cite[Appendix A]{NV}.

\begin{theorem} \label{thm.serreduality}
Suppose that $C$ is a coherent ASF$^{++}$-regular $\ZZ$-algebra of dimension $d \geq 1$ and of Gorenstein parameter $\ell$ with Nakayama isomorphism $\nu:C\to C(-\ell)$.
For $X, Y \in \sD^b(\grmod C)$,
there is a bifunctorial isomorphism
$$
\Hom_{\sD^b(\tails C)}(\pi X, \pi (Y(-\ell))_{\nu}[d-1]) \cong D\Hom_{\sD^b(\tails C)}(\pi Y, \pi X)
$$
i.e. the autoequivalence $S=(-)(-\ell)_{\nu}[d-1]:\sD^b(\tails C) \to \sD^b(\tails C)$ is a Serre functor.
\end{theorem}

\begin{proof}
The fact that the equivalence $(-)(-\ell)_{\nu}:{\GrMod C} \rightarrow {\GrMod C}$ descends to an autoequivalence of $\tails C$ is a straightforward exercise using \cite[Lemma 1.1]{spsmith}.  It follows that there is an induced autoequivalence $S$ on $\sD^b(\tails C)$ as in the statement of the theorem.

First note that, for $X\in \sD^b(\grmod C)$ and for any $n\in \ZZ$, $X_{\geq n}\in \sD^b(\grmod C)$ and $\pi X\cong \pi X_{\geq n}$.  Let $\sD := \sD(\GrMod C)$ and $\sC:=\sD(\Tails C)$.  Since $\operatorname{R}Q \cong \operatorname{R}\omega \circ \pi:\sD\to \sD$ by the proof of \cite[Lemma 4.1.6]{BV} and $(\pi, \operatorname{R}\omega)$ is an adjoint pair of functors between $\sD$ and $\sC$, we have bifunctorial isomorphisms

$$\Hom_{\sC}(\pi X, \pi (Y(-\ell))_{\nu}[d-1]) \cong \Hom_{\sD}(X_{\geq n}, \R Q(Y(-\ell))_{\nu})[d-1]),$$
and  
$$D \operatorname{Hom}_{\sC}(\pi Y, \pi X) \cong D \operatorname{Hom}_{\sD}(Y, \operatorname{R}Q (X_{\geq n}))$$
for any $n\in \ZZ$, so it is enough to show 
$$D \operatorname{Hom}_{\sD}(Y, \operatorname{R}Q (X_{\geq n}))\cong \Hom_{\sD}(X_{\geq n}, \R Q(Y(-\ell))_{\nu})[d-1])$$
for some $n\in \ZZ$ by Lemma \ref{lem.BVl}.  


By Theorem \ref{thm.LD3}, we have a bifunctorial isomorphism
$$D\operatorname{R}\Hom_C(I_{0}(Y), \operatorname{R}QI_{0}(X_{\geq n})) \cong \operatorname{R}\Hom_C(I_{0}(X_{\geq n}), I_{0}(Y) \intotimesc^{\operatorname{L}} D \operatorname{R}Q(C)),$$
so, taking zeroeth cohomology on both sides, 
$$D\Hom_{\sD}(Y, \R Q(X_{\geq n}))\cong \Hom_{\sD}(X_{\geq n}, Y\ulotimes_C D\R Q(C))$$
by Lemma \ref{lem.k-c} (3). 

In order to compute this last expression, we first note that by Lemma \ref{q.ww}, there is a triangle
$$
DC[-1] \rightarrow D\R \t(C)[-1]\cong C(0,-\ell)_{\nu}[d-1] \rightarrow D \operatorname{R}Q(C)
$$
in $\sD(\Bimod (C-C))$. 
Since $(Y \intotimesc^{\operatorname{L}} DC[-1])_{\geq n}=0$ for some $n\in \ZZ$,  we have
$\Hom_{\sD}(X_{\geq n}, Y \intotimesc^{\operatorname{L}} DC[-1])=0$ by Lemma \ref{lem.lln} (1). 
Since $\operatorname{Hom}_{{\sD}}(X_{\geq n},-)$ is cohomological, 
$$\operatorname{Hom}_{{\sD}}(X_{\geq n}, Y \intotimesc^{\operatorname{L}} D \operatorname{R}Q(C)) 
\cong \operatorname{Hom}_{{\sD}}(X_{\geq n}, Y \intotimesc^{\operatorname{L}}C(0,-\ell)_{\nu}[d-1])$$
for any $n\gg 0$.  Finally, since 
\begin{align*}
\operatorname{Hom}_{{\sD}}(X_{\geq n}, Y \intotimesc^{\operatorname{L}}C(0,-\ell)_{\nu}[d-1]) 
& \cong \operatorname{Hom}_{{\sD}}(X_{\geq n}, Y(-\ell)_{\nu}[d-1]) \\
& \cong \Hom_{\sD}(X_{\geq n}, \R Q(Y(-\ell)_{\nu})[d-1])
\end{align*}
for any $n\gg 0$ by Lemma \ref{lem.rell} (3) and Lemma \ref{lem.lln} (2), the result follows. 
\end{proof}

The following lemma is standard.  We include the proof for the convenience of the reader.

\begin{lemma} \label{lem.sgen}
Let $\sC$ be an abelian category and $\cE$ a set of objects in $\sC$.   If $X$ is a bounded complex such that $X^q\in \cE$ for every $q\in \ZZ$, then $X\in \<\cE\>\subset \sD^b(\sC)$.
\end{lemma}

\begin{proof} Since $\<\cE\>$ is closed under shifting complexes, we may assume that $X^q=0$ for every $q<-n$ and $0<q$.  If $n=0$, then the result is trivial.  In general, since there is an exact triangle
$$X^{-n}[n-1]\to X^{\geq -n+1}\to X\to X^{-n}[n],$$
and $X^{-n}[n-1], X^{\geq -n+1}\in \<\cE\>$ by induction, we have $X\in \<\cE\>$.
\end{proof}

\begin{lemma} \label{lem.full}
If $C$ is a right coherent AS-regular $\ZZ$-algebra of dimension $d$ and of Gorenstein parameter $\ell$, then $\<\cP_j, \dots, \cP_{j+\ell-1}\>=\sD^b(\tails C)$ for every $j\in \ZZ$.
\end{lemma}

\begin{proof}
If $C$ is a right coherent AS-regular $\ZZ$-algebra of dimension $d$ and of Gorenstein parameter $\ell$, then
we have exact sequences
\begin{align*}
& 0\to P_{j+\ell} \to F^{d-1}\to \cdots \to F^1 \to P_{j}\to S_{j}\to 0 \\
& 0\to P_{j+\ell-1}\to G^{d-1}\to \cdots \to G^1 \to P_{j-1}\to S_{j-1}\to 0
\end{align*}
in $\grmod C$ where $F^s\in \add\{P_i\}_{j<i<j+\ell}$
for $1\leq s\leq d-1$ and $G^t\in \add \{P_i\}_{j-1<i<j+\ell-1}$
for $1\leq t\leq d-1$ by Remark \ref{rem.lras}, which induce exact sequences
\begin{align*}
& 0\to \cP_{j+\ell} \to \cF^{d-1}\to \cdots \to \cF^1\to \cP_{j}\to 0 \\
& 0\to \cP_{j+\ell-1}\to \cG^{d-1}\to \cdots \to \cG^1\to \cP_{j-1}\to 0
\end{align*}
in $\tails C$, so $\cP_{j-1}, \cP_{j+\ell}\in \<\cP_j, \dots, \cP_{j+\ell-1}\>$ by Lemma \ref{lem.sgen}.
By induction, $\cP_i\in \<\cP_j, \dots, \cP_{j+\ell-1}\>$ for every $i\in \ZZ$.
Since every $X\in \cD^b(\grmod C)$ has a finitely generated free resolution of finite length by Lemma \ref{lem.sgl2},
$$
\<\cP_j, \dots, \cP_{j+\ell-1}\>=\<\{\cP_i\}_{i\in \ZZ}\>=\cD^b(\tails C)
$$
by Lemma \ref{lem.sgen}.
\end{proof}


The theorem below is the converse to Theorem \ref{thm.cs1}, a generalization of the other direction of \cite[Theorem 4.1]{MU1}.  (See also \cite[Theorem 4.7]{Mbc}).

\begin{theorem} \label{thm.cs++2}
If $C$ is a right coherent ASF$^{++}$-regular $\ZZ$-algebra of dimension $d\geq 1$
and of Gorenstein parameter $\ell$ with the Nakayama isomorphism $\nu:C\to C(-\ell)$, then
\begin{enumerate}
\item[(GH1)] $(-)\otimes _{\cC}\omega _{\cC}:=(-)(-\ell)_{\nu}$ is a canonical bimodule for $\tails C$, and
\item[(GH2)] $\{\cP_{-i}\}_{i\in \ZZ}$ is an ample sequence for $\tails C$ which is a full geometric helix of period $\ell$ for $\cD^b(\tails C)$,
\end{enumerate}
so that $\tails C$ satisfies (GH) of period  $\ell$.
\end{theorem}

\begin{proof}
(GH1): This follows from Theorem \ref{thm.serreduality}.

(GH2): By Lemma \ref{lem.aschi}, $\{\cP_{-i}\}_{i\in \ZZ}$ is an ample sequence for $\tails C$.  Since $C$ is $\ell$-periodic, $\cP_{-j-\ell}\cong \cP_{-j}(\ell)_{\nu^{-1}}\cong \cP_{-j}\lotimes _{\cC}\omega _{\cC}^{-1}$ for every $j\in \ZZ$ by Lemma \ref{lem.pper}, so $\{\cP_{-i}\}_{i\in \ZZ}$ satisfies (H1).
Since $(D(C(0,-\ell)))_{ij}=D(C(0,-\ell)_{ji})=D(C_{j, i-\ell})$,  if $d=1$, then
\begin{eqnarray*}
\Ext_{\cC}^q(\cP_j, \cP_i) & \cong & \uExt^q_{\cC}(\cC, \cC)_{ij} \\
& \cong & \begin{cases}
(C\oplus D(C(0,-\ell)))_{ij}=C_{ij}\oplus D(C_{j, i-\ell}) & \textnormal { if $q=0$ } \\
0 & \textnormal { otherwise,} \end{cases}
\end{eqnarray*}
and,  if $d\geq 2$, then
$$\Ext_{\cC}^q(\cP_j, \cP_i)\cong \uExt^q_{\cC}(\cC, \cC)_{ij}\cong \begin{cases} C_{ij} & \textnormal { if $q=0$ } \\
(D(C(0,-\ell)))_{ij}=D(C_{j, i-\ell}) & \textnormal { if $q=d-1$ } \\
0 & \textnormal { otherwise} \end{cases}$$
by
Proposition \ref{prop.asf2'}.
Since $D(C_{i, i-\ell})=0$,  we have
$$
\Ext^q_{\cC}(\cP_i, \cP_i)\cong \begin{cases} C_{ii}=k & \textnormal { if } q=0, \\
0 & \textnormal { if } q\neq 0\end{cases}
$$
for every $i$, so $\{\cP_{-i}\}_{i\in \ZZ}$ satisfies (H2).   If $j<i$, then $C_{ij}=0$ and if $i<j+\ell$, then $D(C_{j, i-\ell})=0$, so, in either case, $\Ext^q_{\cC}(\cP_j, \cP_i)=0$ for every $q$ and $j<i< j+\ell$, so $\{\cP_{-i}\}_{i\in \ZZ}$ satisfies (H3).  Moreover,
$\Ext^q_{\cC}(\cP_j, \cP_i)=0$ for every $q\neq 0$ and every $i\leq j$, so $\{\cP_i\}_{i\in \ZZ}$ is a geometric helix of period $\ell$ for $\cD^b(\tails C)$.
By Lemma \ref{lem.full}, $\{\cP_i\}_{i\in \ZZ}$ is full.
\end{proof}


\subsection{An Application to Noncommutative Quadric Hypersurfaces}

For the rest of the paper, we assume that $k$ is an algebraically closed field of characteristic 0.
If $C$ is a 3-dimensional ``cubic" AS-regular $\ZZ$-algebra $C$ in the sense of \cite{Vq},  which in particular implies that $C$ is a right noetherian AS-regular $\ZZ$-algebra of dimension 3 and of  Gorenstein parameter 4 (cf. \cite [Corollary 5.5.9]{Vq}), then $\tails C$ is considered as a noncommutative $\PP^1\times \PP^1$.  On the other hand, if $S$ is a 4-dimensional noetherian quadratic AS-regular algebra, and $f\in S_2$ is a regular normal element, then $\tails S/(f)$ is considered as a noncommutative quadric surface (see \cite{SV}).  Since $\PP^1\times \PP^1$ is isomorphic to a smooth quadric in $\PP^3$ in commutative algebraic geometry, we expect a similar result in noncommutative algebraic geometry. The following implication was known:

\begin{theorem} \label{thm.a1} \textnormal {\cite [Corollary 6.7]{Vq}}  For every 3-dimensional ``cubic" AS-regular $\ZZ$-algebra $C$, there exist a 4-dimensional right noetherian quadratic AS-regular
algebra $S$ and a regular normal element $f\in S_2$ such that $\tails C\cong \tails S/(f)$.
\end{theorem}

This paper gives a partial converse.

\begin{theorem} \label{thm.a2}
If $S$ is a 4-dimensional
noetherian quadratic AS-regular
algebra and $f\in S_2$ is a regular central element such that
\begin{itemize}
\item{} $S/(f)$ is a domain,
\item{} $S/(f)$ is a noncommutative graded isolated singularity in the sense that $\gldim (\tails S/(f))<\infty$ (that is, $\tails S/(f)$ is ``smooth''),
and
\item{} $S/(f)$ is ``standard" in the sense of \cite[Section 5]{MU1},
\end{itemize}
then
there exists a right noetherian AS-regular $\ZZ$-algebra $C$ of dimension 3 and of Gorenstein parameter 4 such that $\tails S/(f)\cong \tails C$.
\end{theorem}

\begin{proof}
It is known that $\gldim S/(f)=2$ (cf. \cite[Section 5]{MU1}).  By \cite [Theorem 5.16]{MU1}, there exists a full geometric helix $\{E_i\}_{i\in \ZZ}$ of period 4 for $\sD(\tails S/(f))$.  By the proof of \cite[Theorem 5.17]{MU1} and Lemma \ref{lem.aqv}, $\{E_i\}_{i\in \ZZ}$ is ample, so $\tails S/(f)$ satisfies (GH) of period 4, hence there exists an AS-regular $\ZZ$-algebra $C$ of dimension 3 and  of Gorenstein parameter 4 such that $\tails S/(f)\cong \tails C$ by Theorem  \ref{thm.cs1}.   Since $S/(f)$ is right noetherian, $\tails S/(f)$ is a noetherian categoy,  so $C$ is right noetherian
by Theorem \ref{thm.npl}.
\end{proof}

We partially extend the above theorem to noncommutative quadric hypesurfaces below.

\begin{theorem} \label{thm.nqh} Let $S:=k\<x_1, \dots, x_n\>/(x_ix_j-\e_{ij}x_jx_i)$ be a $\pm 1$ skew polynomial algebra where $\e_{ii}=1$ for every $i$, $\e_{ij}=\e_{ji}=\pm 1$ for every  $i\neq j$, and $A:=S/(x_1^2+\dots +x_n^2)$.  If $n\geq 3$, then there exists a right noetherian AS-regular $\ZZ$-algebra $C$ of dimension $n-1$ such that $\tails A\cong \tails C$.  In particular, for every (commutative) smooth quadric hypersurface $Q\subset \PP^{n-1}$, there exists a right noetherian AS-regular $\ZZ$-algebra $C$ of dimension $n-1$ and of Gorenstein parameter $\ell=\begin{cases} n-1 & \textnormal { if $n$ is odd,} \\
n & \textnormal { if $n$ is even,}\end{cases}$ such that $\tails C\cong \coh Q$, the category of coherent sheaves on $Q$.
\end{theorem}


\begin{proof}  Since $A$ is right noetherian, $\tails A$ is a noetherian category,  so it is enough to show that $\tails A$ satisfies (GH) by Theorem \ref{thm.cs1} and Theorem \ref{thm.npl}.

(GH1): It is known that $A$ is a noetherian AS-Gorenstein algebra of dimensiona $n-1\geq 2$ (cf. \cite[Section 2.1]{U}) and of Gorenstein parameter $n-2\geq 1$ and that $\gldim(\tails A)=n-2$ (cf. \cite[Section 2.2]{U}),
so $\tails A$ has a canonical bimodule $\omega _{\cA}=\cA_{\nu}(-n+2)$ where $\nu$ is the Nakayama automorphism of $A$ (cf. \cite[Lemma 2.2]{U}.)


(GH2): Let $\operatorname {Ind} ^0(\CM^{\ZZ}(A))=\{A, X_1, X_2, \dots, X_{\a}\}$ be the set of complete representatives of isomorphism classes of indecomposable graded maximal Cohen-Macaulay right $A$-modules generated in degree 0 where we say that $M\in \grmod A$ is a graded maximal Cohen-Macaulay if $\uExt^q_A(M, A)=0$ for every $q\neq 0$.
We label the sequence
$$\cA(-n+3), \dots, \cA(-1), \cA, \cX_1, \cX_2, \dots, \cX_{\a}$$
by $E_0, \dots, E_{\ell-1}$ where $\ell=n-2+\a$.
We extend in both directions the sequence $E_0, \dots, E_{\ell-1}$ to $\{E_i\}_{i\in \ZZ}$ by $E_{i+r\ell}:=E_i\lotimes _{\cA}(\omega _{\cA}^{-1})^{\otimes r}$ so that $\{E_i\}_{i\in \ZZ}$ satisfies (H1).
By \cite [Lemma 3.15]{U} and Lemma \ref{lem.aqv},  $\{E_i\}_{i\in \ZZ}$ is an ample sequence for $\tails A$.


If $X\in \operatorname {Ind} ^0(\CM^{\ZZ}(A))$, then
$X\otimes _A\omega _A^{\otimes r}(r(n-2))\cong X_{\nu^r}\in \operatorname {Ind} ^0(\CM^{\ZZ}(A))$ for every $r\in \ZZ$,
so,  for every $i\in \ZZ$, there exists $s\in \ZZ$ and $X\in \operatorname {Ind} ^0(\CM^{\ZZ}(A))$ such that  $E_i\cong \cX(s)$.
Since $\cA(i)\lotimes _{\cA}\omega _{\cA}^{\otimes r}\cong \cA_{\nu^r}(i-r(n-2))\cong \cA(i-r(n-2))$ for every $i\in \ZZ$, there exists a permutation $\s$ on $\{1, \dots, \a\}$ such that $\cX_j\lotimes _{\cA}\omega _{\cA}^{\otimes r}\cong \cX_{\s^r(j)}(-r(n-2))$ (cf. the proof of \cite[Lemma 3.15]{U}).
It follows that the sequence $\{E_i\}_{i\in \ZZ}$ looks like
\begin{small}
\begin{align*}
\cdots,  \;
& \cA(r(n-2)-n+3), \dots, \cA(r(n-2)-1), \cA(r(n-2)),  \\
& \cX_{\s^{-r}(1)}(r(n-2)), \cX_{\s^{-r}(2)}(r(n-2)), \dots, \cX_{\s^{-r}(\a)}(r(n-2)),  \\
&  \cA((r+1)(n-2)-n+3)=\cA(r(n-2)+1), \dots, \cA((r+1)(n-2)-1), \cA((r+1)(n-2)),  \\
& \cX_{\s^{-(r+1)}(1)}((r+1)(n-2)), \cX_{\s^{-(r+1)}(2)}((r+1)(n-2)), \dots, \cX_{\s^{-(r+1)}(\a)}((r+1)(n-2)), \;
\cdots
\end{align*}
\end{small}

Since $E_0, \dots, E_{\ell-1}$ is an exceptional sequence by \cite[Lemma 3.12]{U},
$$
\Ext^q_{\cA}(E_i, E_i)\cong \begin{cases} k & \textnormal { if } q=0, \\
0 & \textnormal { if } q\neq 0, \end{cases}
$$
so $\{E_i\}_{i\in \ZZ}$ satisfies (H2).


Since $\gldim (\tails A)=n-2$, $\Ext^q_{\cA}=0$ for every $q\geq n-1$.
For $A\not \cong X, Y\in \operatorname {Ind} ^0(\CM^{\ZZ}(A))$,
$$\begin{array}{llll}
(1) & \Ext^q_{\cA}(\cA(s), \cA(t)) & \cong \begin{cases} A_{t-s} & \textnormal { if } q=0, \\
0 & \textnormal { if } 1\leq q\leq n-3, \\
D(A_{s-t-n+2}) & \textnormal { if } q=n-2, \end{cases} \\
(2) & \Ext^q_{\cA}(\cA(s), \cX(t)) & \cong \begin{cases} X_{t-s} & \textnormal { if } q=0, \\
0 & \textnormal { if } 1\leq q\leq n-3, \\
D(\uHom_A(X, A)_{s-t-n+2}) & \textnormal { if } q=n-2, \end{cases} \\
(3) & \Ext^q_{\cA}(\cX(s), \cA(t)) & \cong \begin{cases} \uHom_A(X, A)_{t-s} & \textnormal { if } q=0, \\
0
& \textnormal { if } 0\leq q\leq n-3, \\
D((X_{\nu})_{s-t-n+2}) & \textnormal { if } q=n-2, \end{cases} \\
(4) & \Ext^q_{\cA}(\cX(s), \cY(t)) & \cong \begin{cases} \uExt^q_{A}(X, Y)_{t-s} & \textnormal { if } 0\leq q\leq n-3, \\
D(\uHom_A(Y, X_{\nu})_{s-t-n+2}) & \textnormal { if } q=n-2, \end{cases}
\end{array}$$
by \cite[Lemma 2.3]{U} (cf. the proof of \cite[Lemma 3.13]{U}).
We also use the following facts:
\begin{enumerate}
\item[(i)]
If $i\leq 0$, then $\uHom_A(X, A)_i=0$ (\cite [Lemma 3.7]{U}).
\item[(ii)]
If $i<0$, or $i\leq 0$ and $X\ncong Y$, then $\uHom_A(X, Y)_i=0$ (\cite [Lemma 3.9]{U}).
\item[(iii)] If $q\geq 1$ and $i\neq -q$, then $\uExt^q_A(X, Y)_i=0$ (\cite [Lemma 3.8]{U}).
\end{enumerate}
We will now show that $\Ext^q_{\sC}(E_i, E_j)=0$ for every $q$ when $0<i-j<\ell$.
It is enough to consider the following cases where $A\not \cong X, Y\in \operatorname {Ind} ^0(\CM^{\ZZ}(A))$.
\begin{enumerate}
\item{}  The case $E_i=\cA(s), E_j=\cA(t)$ for some $s, t\in \ZZ$:  If $0<i-j<\ell$, then $0<s-t<n-2$. Since $t-s<0$ and $s-t-n+2<0$, we have $\Ext^q_{\cA}(\cA(s), \cA(t))=0$ for every $q\in \ZZ$.
\item{}  The case $E_i=\cA(s), E_j=\cX(t)$ for some $s, t\in \ZZ$: If $j<i<j+\ell$, then $E_i=\cA(s), E_j=\cX(t)$ are positioned as follows
\begin{align*}
& \cX_{\s^{-r}(1)}(r(n-2)), \dots \cX(t)=\cX(r(n-2)), \dots, \cX_{\s^{-r}(\a)}(r(n-2)),  \\
& \cA(r(n-2)+1), \dots, \cA(s), \dots,  \cA((r+1)(n-2)),
\end{align*}
that is, $t=r(n-2)$ and $r(n-2)+1\leq s\leq (r+1)(n-2)$ for some $r\in \ZZ$, so $t<s\leq t+n-2$.
Since $t-s<0$ and $s-t-n+2\leq 0$, we have $\Ext^q_{\cA}(\cA(s), \cX(t))=0$ for every $q\in \ZZ$ by (i).
\item{} The case $E_i=\cX(s), E_j=\cA(t)$ for some $s, t\in \ZZ$: If $j<i<j+\ell$, then $E_i=\cX(s), E_j=\cA(t)$ are positioned as follows
\begin{align*}
& \cA(r(n-2)-n+3), \dots, \cA(t), \dots, \cA(r(n-2)),  \\
& \cX_{\s^{-r}(1)}(r(n-2)), \dots \cX(s)=\cX(r(n-2)), \dots, \cX_{\s^{-r}(\a)}(r(n-2)),
\end{align*}
that is, $r(n-2)-n+3\leq t\leq r(n-2)$ and $s= r(n-2)$ for some $r\in \ZZ$, so $t\leq s< t+n-2$.
Since $t-s\leq 0$ and $s-t-n+2<0$, we have  $\Ext^q_{\cA}(\cX(s), \cA(t))=0$ for every $q\in \ZZ$ by (i).
\item{} The case $E_i=\cX(s), E_j=\cY(t)$ for some $s, t\in \ZZ$:
If $j<i<j+\ell$, then $E_i=\cX(s), E_j=\cY(t)$ are positioned as either
\begin{small}
$$\cX_{\s^{-r}(1)}(r(n-2)), \dots, \cY(t)=\cY(r(n-2)), \dots,  \cX(s)=\cX(r(n-2)), \dots, \cX_{\s^{-r}(\a)}(r(n-2))$$
\end{small}
or
\begin{small}
\begin{align*}
& \cX_{\s^{-r}(1)}(r(n-2)), \dots, \cX_{\nu}(r(n-2)), \cdots, \cY(t)=\cY(r(n-2)), \dots, \cX_{\s^{-r}(\a)}(r(n-2)),  \\
& \cA(r(n-2)+1), \dots, \cA((r+1)(n-2)-1), \cA((r+1)(n-2)),  \\
& \cX_{\s^{-(r+1)}(1)}((r+1)(n-2)), \dots, \cX(s)=\cX((r+1)(n-2)), \dots, \cX_{\s^{-(r+1)}(\a)}((r+1)(n-2)),
\end{align*}
\end{small}
so either $s=t$ and $Y\not \cong X$, or $s=t+n-2$ and $Y\ncong X\otimes _A\omega_A(n-2)\cong X_{\nu}$.

(a) The case $s=t$ and $Y\ncong X$:
If $0\leq q\leq n-3$, then
$$
\Ext^q_{\cA}(\cX(s), \cY(t))\cong \uExt^q_A(X, Y)_0=0
$$
by (ii) and (iii) since $Y\ncong X$.  If $q=n-2$, then $\Ext^{n-2}_{\cA}(\cX(s), \cY(t))\cong D(\uHom_A(Y, X_{\nu})_{s-t-n+2})=0$ by (ii) since $s-t-n+2=-n+2<0$.

(b) The case $s=t+n-2$ and $Y\ncong
X_{\nu}$:
If $0\leq q\leq n-3$, then $\Ext^q_{\cA}(\cX(s), \cY(t))\cong \uExt^q_A(X, Y)_{t-s}=0$ by (ii) and (iii) since $t-s=-n+2<0$ and $t-s=-n+2\neq -q$.  If $q=n-2$, then $\Ext^{n-2}_{\cA}(\cX(s), \cY(t))\cong D(\uHom_A(Y, X_{\nu})_0)=0$ by (ii) since $Y\ncong X_{\nu}$.
\end{enumerate}
It follows that $\{E_i\}_{i\in \ZZ}$ satisfies (H3).


On the other hand, if $i\leq j$ and $E_i=\cX(s), E_j=\cY(t)$ for some $X, Y\in \operatorname {Ind} ^0(\CM^{\ZZ}(A))$ and $s, t\in \ZZ$, then $s\leq t$, so $t-s\geq 0$ and $s-t-n+2<0$, hence
$\Ext^q_{\cA}(\cX(s), \cY(t))=0$ for every $q\neq 0$  by (i), (ii), and (iii).
It follows that
$\{E_i\}_{i\in \ZZ}$ is a geometric helix of period $\ell$ for $\sD^b(\tails A)$.

Since $E_0, \dots, E_{\ell-1}$ is full by \cite[Lemma 3.12]{U}, $\{E_i\}_{i\in \ZZ}$ is full by \cite[Lemma 3.16, Remark 3.17]{MU1}.

It follows that $C:=C(\tails A, \{E_i\}_{i\in \ZZ})$ is a right noetherian AS-regular $\ZZ$-algebra of dimension $n-1$ such that $\tails A\cong \tails C$ by Theorem \ref{thm.cs1} and Theorem \ref{thm.npl}.

If $Q\subset \PP^{n-1}$ is a (commutative) smooth quadric hypersurface, then
$$
Q\cong \Proj k[x_1, \dots, x_n]/(x_1^2+\cdots +x_n^2),
$$
so
$$
\coh Q\cong \tails k[x_1, \dots, x_n]/(x_1^2+\cdots +x_n^2),
$$
hence the final assertion (cf. \cite[Theorem 1.1]{U}).
\end{proof}

\begin{remark}
\begin{enumerate}
\item{} Presumably, the AS-regular $\ZZ$-algebra $C$ constructed above is not 1-periodic
(see the quiver presentations of $\End_{\cA}(\oplus _{i=0}^{\ell-1}E_i)$ in \cite [Section 3.5]{U}), so that $C$ is not a $\ZZ$-algebra associated to any (AS-regular) graded algebra by Lemma \ref{lem.1p}.
\item{} Let $Q\subset \PP^{n-1}$ be a (commutative) smooth quadric hypersurface.  If $n$ is odd, then the final assertion of the above theorem follows from \cite[Proposition 3.3]{BP}.   If $n$ is even, then the $\ZZ$-algebra $C$ constructed in the above theorem is AS-regular, but not Koszul.
\end{enumerate}
\end{remark}


\end{document}